%
%  VERSION  28. November 1999
%
%  geaenderte outputparameter fuer doppelseitiges
%  kopieren
%
%
\widowpenalty=15000
\overfullrule=0pt
\pretolerance=-1
\tolerance=2500
\doublehyphendemerits 50000
\finalhyphendemerits 25000
\adjdemerits 50000
\hbadness 1500
\abovedisplayskip=7pt plus 5pt minus 6pt
\belowdisplayskip=7pt plus 5pt minus 6pt
\hsize39pc \vsize52pc
\baselineskip11.5pt
\magnification=1200
\parindent 0pt
\def\exp{{\rm \hskip1.5pt  exp \hskip1.5pt }}
\def\summe{\sum\nolimits}
\font\gross=cmbx12
\font\grossrm=cmr12
\font\gr=cmr17

%\font\nineit=amti9
%\font\eightit=amti8

%
%\font\calligX=callig15 at 10pt
%\newfam\calligfam
%\textfont\calligfam=\calligX
%\font\calligVIII=callig15 at 8pt
%\scriptfont\calligfam=\calligVIII
%\font\calligV=callig15 at 5pt
%\scriptscriptfont\calligfam=\calligV
%\def\script{\fam\calligfam\calligX}
%
\def\parein#1#2{\par\noindent\rlap{\rm {#2}}\parindent#1\hang\indent\ignorespaces}
\def\paraus{\parindent 0pt\par}

\def\G{\Gamma}

\def\a{\alpha}

\def\l{\ell}
\def \ol{\overline}
\def \1{\backslash}
%\input NEWMAKROS.tex
%
%
% output parameter
%
\output={
\ifodd\pageno\hoffset=0pt
\else
\hoffset=0truecm\fi\plainoutput}
\font\authorfont=cmcsc10 at 12pt
\font\titlefont=ptmb at 13pt%cmbx12 at 14pt 
\font\affiliationfont=cmti12 at 11pt
 
\font\gross=cmbx12

\font\ninecaps=cmcsc9
\font\ninebf=cmbx9
\font\nineit=cmti9
\font\ninerm=cmr9
\font\eightrm=cmr8
\font\eightit=cmti8
\font\eightbf=cmbx8
\def\er{{\rm Erfc}\,}
\def\re{{\rm Re}\,}
\def\im{{\rm Im}\,}
\def\la{{\bf L}}
\def\fa{{\bf F}}
\def\ia{{\bf I}}
\def\ctrig{ C_{\rm trig}}
\def\chyp #1 { C_{{\rm hyp}, #1 }}

\abovedisplayskip=8pt plus 5pt minus 3pt
\belowdisplayskip=8pt plus 5pt minus 3pt
\baselineskip=15pt
%
%\centerline{  } \pageno=0 \vfill\eject
%
\centerline{\titlefont ON THE VALUATION OF ARITHMETIC-AVERAGE}
\centerline{\titlefont ASIAN OPTIONS: INTEGRAL REPRESENTATIONS}
\baselineskip=11.5pt
\vskip.3cm
\centerline{\authorfont Michael  Schr\"oder}
\vskip.2cm
\centerline{\affiliationfont Lehrstuhl Mathematik III}
%\vskip.05cm
\centerline{\affiliationfont Seminargeb\"aude A5,
Universit\"at Mannheim,  D--68131 Mannheim}
\vskip.5cm 
\centerline{
\vbox{
\hsize=9cm
\baselineskip=9truept
\ninerm
This paper derives  integral representations for the Black--Scholes 
price of arithmetic--average Asian options.  Their proof is by Laplace 
inverting  the Laplace transform of [{\ninebf GY}] using complex analytic 
methods. The analysis  ultimately rests on the gamma function which 
in this sense is at the base of Asian options. The results of 
[{\ninebf GY}] are corrected and their validitity is extended.}
}
\vskip.5cm
{\bf Introduction:}\quad The aim of this paper is to derive exact 
closed--form valuation formulas for European--style arithmetic--average Asian 
options in a Black--Scholes framework. These options are a particular class of 
path--dependent options with the arithmetic average of the underlying
security as the relevant variable. 
\medskip
This paper has its origin in that developement in the analysis of Asian 
options initiated by~[{\bf Y}]. Yor's valuation formula 
gives clear evidence that pricing Asian options 
is a problem of some intrinsic difficulty indeed for which no, in the strict 
sense, simple solution should be expected. Instead, one should, 
as a first step, ask for structurally clear solutions, and only then, 
as a second step, consider in particular questions of a more computational 
nature.     
\medskip
The main  valua\-tion result of this paper is described in paragraph four 
and identifies the Black-Scholes price of an Asian option as  
an integral over, roughly speaking, the product of two well--studied 
higher transcendental functions both given as integrals 
and built up using so--called Hermite functions. These last functions   
generalize the familiar 
complementary error function and come from boundary--value 
problems in potential theory for domains whose surface is an infinite 
parabolic cylinder. It seems to be the first time that such parabolic cylinder 
functions are explicitly identified to characterize solutions to 
problems of Finance. Their relevant properties are reviewed in paragraph 
zero. 
\medskip
Both valuation formulas involve three integrations. In a sense, this 
measures the difficulty of the problem. For getting some intuition, first 
think of the error integral as an example of how a single integration 
yields a higher transcendental function. In the best case, a triple integral 
separates as a product of three such functions. In the worst case, no  
separation is possible and all one can say is that one has to integrate a 
function in three variables three times in succession. 
With its integrand a function $f(x,y,w)$ that cannot be further separated, 
Yor's formula belongs to this last category. In our formula complete 
separation is not achieved either. However, with its integrand a product  
$g(x,y)h(y,w)$ it is a single integral of 
a function given as the product of two single integrals. 
\medskip
The method used in this paper is the Laplace transform Ansatz to option 
pricing  of [{\bf GY}] and its main result is proved by Laplace inverting 
the Laplace transform they have derived.  
\medskip
Transform methods like the Laplace transform are classical concepts
in analysis with applications ranging from number theory to
boundary--value problems for partial differential equations. 
In particular the Laplace transform approach used in this paper generally 
consists of two stages. As a first step, a suitably nice function on the 
positive real line is transformed into a complex--valued function. This 
Laplace transform exists on a 
half--plane sufficiently deep within the right complex half--plane  
and defines there a holomorphic, i.e., complex analytic, function.
The methods in this part of the transform analysis have real 
analysis origin. 
As a second step, the Laplace transform has 
to be inverted to give the desired function on the positive real line. The 
natural methods in this step are complex 
analytic. Some relevant concepts are recalled in paragraph zero
of this paper
%\goodbreak
\medskip
The computation of the Laplace transform already is quite a problem 
in general. Indeed, with the functions to be transformed generally not 
known one expressed aim for trying to calculate their Laplace 
transform actually is to get an explicit expression at all. However, 
for an actual calculation of a Laplace transform in such a situation 
additional insights into the mathematical structure are necessary. Thus, 
the Laplace transform  in [{\bf GY}] is computed introducing a stochastic 
clock and using general results on Bessel processes. 
\medskip
It should be said that, unfortunately, there are serious problems of 
a fundamental nature with this result. More precisely, the conceptual 
error I found in September 1999 kills this approach to the valuation 
of Asian options as it stands. Luckily, Peter Carr eventually succeeded 
in convincing me that the Geman--Yor result can nevertheless be used in 
a very astute way for its original purpose. All this is discussed in 
paragraph five of this paper. 
\medskip
Mathematically, the inversion of a non--trivial Laplace transform 
belongs to the nastier problems in analysis. Difficulties are encountered 
here as a rule, and more than often one has to be satisfied with only 
knowing the Laplace tranform. It thus can be regarded as as the main
mathematical contribution of the paper that it indeed provides such an 
inversion for the above Laplace transform of the price of Asian options. 
The details of the proof are described in paragraphs six to twelve of 
this paper. Its methods are complex analytic. It uses the complex inversion 
formula for the Laplace transform and is ultimately based on the classical 
Hankel formulas for the gamma function. These methods are such that the 
results can be extended from the classical Black--Scholes framework 
essential for the validity of [{\bf GY}] to one with no restrictions 
on the drift coefficients.   
\medskip
It came as a surprise to the author that one has such a structural integral 
representation for the value of Asian options as derived in this paper. 
He is not sure if this should be attributed to the qualities intrinsic to
Asian options or to the fact that, as described above, one takes a somewhat 
indirect approach to their valuation. 
\nopagenumbers
\def\Verfasser{{\ninecaps Michael Schr\"oder}}
\def\Kapiteltitel{{\ninecaps Valuation of arithmetic--average options}}
\headline={\vbox{\line{\ninebf \ifodd\pageno\hss \Kapiteltitel\hss
\folio
\else
\folio\hss\Verfasser\hss\fi}
%\hrule
}\hss}
\medskip
The series representations and asymptotic expansions of our formula  
appropriate and adequate for computation are discussed in [{\bf Sch}].
\bigskip
\vbox{\baselineskip=8.5pt
\eightrm 
{\eightbf Acknowledgements:}\quad This paper owes its existence very much to
the open and stimulating atmo\-sphere at the number theory chair of Professor
Weissauer. It is a pleasure to thank in par\-ti\-cular R. Weissauer, 
J. Ballmann, and D. Fulea for their support and the genuine interest they
have taken in this project.  I am very grateful to Professor Harder for
making possible my very pleasant stay at the {\eightit Mathematisches 
Institut der Universit\"at Bonn\/}  in the {\eightit Wintersemester 97/98\/}
during which this paper was completed. Meanwhile I also  wish to thank 
Professor Yor for correspondence and his comments. To Dr. P. Carr 
(Banc of America Securities, New York) I owe a special debt for his help 
in salvaging the valuation results of [{\eightbf GY}] in general 
and this paper in  particular.  According to P. Carr I owe in this way 
a debt also to Professor D. Madan (University of Maryland). 
Finally, I wish  to thank very much Professor Pliska for his efforts, 
and for his work as an editor of a paper that proved to be difficult 
to handle. 
%
% 
%$$\vcenter\vbox{\baselineskip=10pt
%\hbox to 11.7cm{\hfil \ninebf Table of Contents\hfil}
%\vskip3pt
%\hrule
%\vskip6pt
%\halign{\hfil{\ninerm #}\quad&{\ninerm #}
%\dotfill&\quad \hfil{\ninerm #}\hfil\cr  
%1.&Black--Scholes framework for valuating contingent claims&2\cr 
%2.&The notion of arithmetic--average Asian options&2\cr 
%3.&The valuation problem for arithmetic--average Asian options&3\cr 
%4.&Statement of results&3\cr 
%5.&Characterizing the price of Asian options by its Laplace transform
%&5\cr 
%6.&Preliminaries on integral representations of Gamma and Bessel functions$
%\ldots$&7\cr 
%7.&First steps of the Laplace inversion&9\cr 
%8.&Calculation of certain Laplace transforms&11\cr 
%9.&Two intermediate results&13\cr
%10.&Proof of the valuation formula&14\cr 
%11.&References&16\cr 
%}
%}}$$  
}
\goodbreak
\bigskip
{\bf 0.\quad Analytic preliminaries:}\quad This paragraph collects some 
relevant facts about the Laplace transform from [{\bf B}], [{\bf D}], and 
about Hermite functions from [{\bf L}, \S \S 10.2ff]. 
\bigskip
{\bf Laplace transform:}\quad 
The class of functions considered in the sequel is that of 
functions of {\it exponential type\/}, i.e., of continuous, real--valued
functions $f$ on the non--negative real line 
such that there is a real number $a$ for which $\exp(at)f(t)$ is 
bounded for any $t>0$. The {\it Laplace transform\/} is the  linear 
operator $\la$ that associates to any function $f$ of exponential type 
the complex--valued  function  $\la(f)$, holomorphic on 
a suitable complex half--plane, which for any complex number 
$z$ with sufficiently big real part is explicitly given by:
$$ \la(f)(z)=\int\nolimits_0^\infty e^{-zt}f(t)\,dt.
$$
\goodbreak
This operator is an injection. 
Its inverse, the inverse Laplace transform  $\la ^{-1}$, is expressed 
as a contour integral by the {\it complex inversion formula\/} of 
Riemann. Applying to any function $H$ analytic on 
half--planes $\{ \re (z)\!\ge\!z_0\}$ with $z_0$ any sufficiently 
big positive real number, it asserts: 
$$ \la^{-1}(H)(t)={1\over 2\pi i}\int\nolimits_{z_0-i\infty}^{z_0+i\infty}
e^{zt}\cdot H(z)\, dz,$$ 
for any positive real number $t$ if $H$ satisfies a  growth 
condition at infinity such that the above integral exists. For any 
function $f$ of exponential type one so has in particular %: 
$(\la ^{-1}\circ \la) (f)=f$.  
\medskip
For getting an idea of the proof of the complex inversion formula 
start with any function $f$ holomorphic on a 
half--plane $\{ z| \re (z)>a\}$. For obtaining a suitable 
integral representation of $f$, invoke the Cauchy integral formula 
in the following way. Fix $z$ with $\re (z)>a$ and choose a line $C(z_0)=
\{ z| \re (z)=z_0\}$ such that $z$ is to the right of it, i.e., $z_0<\re(z)$,
and such that it is contained in the half--plane where $f$ is holomorphic, 
i.e. $a<z_0$. Make this line into a path of integration by moving upwards 
from $z_0-i\infty$ to $z_0+i\infty$. Shift pespective to the Riemann sphere by 
adding the point infinity to the complex plane. In this picture, $C(z_0)$
is a closed path that circles clockwise, i.e., in the mathematically 
negative direction, around $z$. Formally apply the Cauchy integral formula,
expressing $2\pi i\cdot f(z)$ as minus the integral of $(w-z)^{-1} f(w)$
over $C(z_0)$. Since 
%$ f(z)=-(2\pi i)^{-1} \int_{C(z_0)}(w-z)^{-1} f(w)\, dw$.
%$$ f(z)=-{1\over 2\pi i} \int_{C(z_0)} {f(w)\over w-z}\, dw.$$
%Notice that $-(w-z)^{-1}$ is the Laplace transform at the point $z$ of 
%the map on the positive real line given by $\exp(wh)$ for any $h>0$, whence 
$-(w-z)^{-1}$ is the Laplace transform at $z$ of 
the map on the positive real line given by $\exp(wh)$ for any $h>0$, 
it so follows:
$$ f(z)={1\over 2\pi i} \int_{C(z_0)}\int_0^\infty 
e^{(w-z)h} f(w)\, dh\, dw.$$
On the Riemann sphere
the contour $C(z_0)$ can be seen as the limit of closed contours contained 
in the complex plane. The above integral thus makes sense as such a limit
if $f$ is supposed rapidly decreasing, i.e., the absolute value of $f$ 
exponentially decreases to zero with the absolute value $|w|$ of $w$ going 
to infinity. Then apply Fubini's theorem to interchange the integrals in the 
above formula.
%To make sense of this integral suppose $f$ is rapidly decreasing, i.e., the 
%absolute value of $f$ exponentially decreases to zero with the absolute
%value $|w|$ of $w$ going to infinity. Apply Fubini's theorem 
%to interchange the integrals in the 
%above formula.
%%$$ f(z)=\int_0^\infty e^{-zh}\,  {1\over 2\pi i} \int{C(z_0)}
%%e^{wh} f(w)\, dw\, dh\, , $$
%%identifies 
The inverse Laplace transform of $f$ is then identified as the function on the 
positive real line given by
$${1\over 2\pi i} \int_{C(z_0)}
e^{wh} f(w)\, dw\, ,$$
for any $h>0$, which is the complex inversion formula. 
%[{\bf B}, Chapter 7, \S 5]. 
\goodbreak
\bigskip
{\bf Hermite functions:}\quad The {\it Hermite function\/} $H_\nu$ 
of degree any complex number $\nu$ is the complex--valued function on the 
complex plane given as solution to the differential equation 
$u''-2zu'+2\nu u=0$ in the complex variable $z$. It is holomorphic
on the complex plane as a function of both its variable $z$ and its 
degree $\nu$. Using the differential equation, 
for $\nu$ any non--negative integer, $H_\nu$ is the $\nu$--th  
Hermite polynomial and belongs to a standard class of orthogonal
functions. If the real part of the degree $\nu$ is negative the Hermite
function $H_\nu$ has the following integral representation:
$$H_\nu(z)={1\over \Gamma(-\nu)}\int_0^\infty e^{-u^2-2zu}u^{-(\nu+1)}\, du \,
 ,$$
for any complex number $z$ and with $\Gamma(-\nu)$ the value of the gamma 
function at $-\nu$. Hermite functions thus also generalize the {\it 
complementary error function\/} $\er$. Recalling that $\er$ is for any 
complex number $z$ given by any of the following expressions:
$$ \er (z)=N(-z)= 
{2\over \sqrt{\pi}} \int _z^{\infty} e^{-u^2}du
=e^{-z^2}{2\over \sqrt{\pi}}\int_0^\infty e^{-u^2-2zu}du
\, ,$$ 
notice $(2/\sqrt{\pi})\cdot H_{-1}(z)=\exp(z^2)\er (z)$. More generally, one 
has for any non--negative integer $n$ the identity:
$$ H_{-(n+1)}(z)= {\sqrt{\pi}\over 2}\cdot {(-1)^n\over n!\cdot 2^n}\cdot 
{d^n\over dz^n}\Bigl( e^{z^2}\er(z)\Bigr)$$
for any complex number $z$. Hermite
functions of degree any complex number $\nu$ satisfy the recurrence relations:
$$\leqalignno{
&H_\nu'(z)=2\nu\cdot H_{\nu-1}(z),&\cr
\noalign{\vskip3pt}
&H_{\nu+1}(z)-2z\cdot H_{\nu}(z)+2\nu\cdot H_{\nu-1}(z)=0.&\cr 
}$$
Using the differential equation of the Hermite functions the second of 
these is implied by the first which can be seen using
the following absolutely convergent series expansion: 
$$H_\nu(z)= {1\over 2\cdot \Gamma(-\nu)} \summe_{n=0}^\infty
{(-1)^n\over n!}\, \Gamma\left( {n-\nu\over 2}\right) (2z)^n, 
$$
valid for any complex number $z$ and uniformly convergent on compact sets.
Like the exponential series this series has excellent convergence properties
for complex numbers $z$ of smaller absolute value. For arguments $z$ of 
larger absolute value one uses  asymptotic expansions. 
If the real part of $\nu $ is negative the Hermite function  
of degree $\nu$ has the following asymptotic expansion in the right 
half--plane:
$$ H_\nu(z)=(2z)^\nu\summe_{k=0}^n (-1)^k{ (-\nu)_{2k}\over k!}\cdot 
{1\over (2z)^{2k}} +r_{n+1}(z),$$
for any complex number $z$ with positive real part and any 
non--negative integer $n$. Herein $(-\nu)_{2k}$ is the Pochhammer symbol 
recalled to be given by $(-\nu)_0=1$ and $(-\nu)_{2k}
=\Gamma(-\nu+2k)/\Gamma(-\nu)=(-\nu)(-\nu+1)\cdots (-\nu+2k-1)$ for $k$
positive. For estimating the error term $r_{n+1}$, fix any positive real 
number $\delta$ less than or equal to $\pi/2$ and choose any complex 
number $z$ in the wedge contained in the right half--plane that is enclosed 
between the two rays emanating from the origin with angles $\pm(\pi/2-\delta)$ 
respectively. Equivalently, the absolute value of the argument of $z$ is 
at most $\pi/2-\delta$. Then the absolute value of the error term 
$r_{n+1}(z)$ satisfies the following estimate:
$$ |r_{n+1}(z)|\le  
{1\over \bigl(\sin(\delta)\bigr)^{2(n+1)-\re (\nu)}}\cdot
{\Gamma\bigl(2(n+1)-\re (\nu)\bigr)\over (n+1)!\, |\Gamma(-\nu)| }\cdot
{1\over (2|z|)^{2(n+1)-\re (\nu)} } \,. $$
The differential equation for the Hermite function
of degree any complex number $\nu$ can be transformed into the differential
equation $ u''+(\nu + 1/2 -z^2/4)u=0$ whose solution $D_\nu$ is thus 
connected with the Hermite function $H_\nu$ by 
$D_\nu(z)= 2^{-\nu/2} \exp({-z^2/4}) H_\nu(2^{-1/2} z)$, for any 
complex number 
$z$. Similarly, Hermite functions are related to the second Kummer
confluent hypergeometric function $\Psi$ by 
$H_\nu(z)=2\cdot\Psi(-\nu/2, 1/2, z^2)$, for any complex number $z$ in
the right half--plane. 
%\goodbreak
\vskip.8cm
\centerline{\gross %PART I \qquad STATEMENT OF RESULTS}
Part I\qquad Statement of results}
\vskip.5cm
{\bf 1.\quad Black--Scholes framework for valuating contingent claims:}\quad
For our analysis we place ourselves in the Black--Scholes framework and use
the risk--neutral approach to valuating contingent claims as described in 
[{\bf DSM}, Chapters 17, 22]. 
\bigskip
In this set--up one has two securities. First there is a riskless security,
a bond, that has the continuously compounding positive interest rate $r$.
Then there is a risky security whose price process $S$ is modelled as 
follows. Start with a complete probability space equipped with the standard 
filtration of a standard Brownian motion on it that has the time set 
$[0,\infty)$. On this filtered space one has the 
{\it risk neutral measure\/} $Q$,
a probability measure equivalent to the given one, and a standard $Q$--Brownian
motion $B$ such that $S$ is the strong solution of the following stochastic 
differential equation:
%
%More precisely, start with a complete probability space 
%$(\Omega, {\cal F}, P)$
%equipped with the standard filtration ${\bf F}$ of a standard Brownian motion
%with time set $[0, \infty)$. Assume there is a riskless security, a bond, 
%that has the continuously compounding non--negative interest rate $r$, i.e., 
%whose price process $\beta$ solves the following differential equation:
%$$ d\beta_t=r\cdot \beta_t\, dt, \qquad t\in [0,\infty).$$
%Moreover assume that there exists a risky security whose price process $S$
%is given by the following stochastic differential equation:
%
$$ dS_t=\varpi\cdot S_t\, dt +\sigma\cdot S_t\, dB_t\, , 
\qquad  t\in [0,\infty).$$ 
%Here $B$ is a standard Brownian motion on the filtered probability space 
%$(\Omega, {\cal F}, {\bf F}, Q)$ with $Q$ a {\it risk--neutral\/} probability
%measure equivalent to $P$ as described in [{\bf DSM}, 17E]. 
Herein the positive constant $\sigma$ is the volatility of $S$, whereas the 
specific form of the constant $\varpi$ depends on the security modelled. For
instance it is the interest rate rate if $S$ is a non--dividend paying stock.
\bigskip
With the existence of $Q$ one has the {\it arbitrage--pricing principle\/}: 
let $Y$  be any European--style contingent claim on the above filtered
probability space written at time $t_0$ and paying $Y_T$ at time $T$. 
At any time $t$ between $t_0$ and $T$ its time--$t$ price $C_t$ is then: 
$$ C_t
%=E^Q\left[ e^{-\int _t^T r \, d\tau}\, Y_T\,\Big\vert\, {\cal F}_t\right]
= e^{-r(T-t)} E_t^Q\bigl[Y_T\bigr]$$
with the expectation $E_t^Q$ conditional on the information 
%${\cal F}_t$ 
available at time $t$ and taken with respect to $Q$.
\bigskip
{\bf 2.\quad The notion of arithmetic--average Asian options:}\quad Fix any 
time $t_0$ and consider in the Black--Scholes framework reviewed above 
the process $J$ given for any time $t$ by:
$$ J(t)=\int_{t_0}^t S_\tau \, d\tau . $$
The European style {\it arithmetic--average Asian option\/} written at
time $t_0$, with maturity the time $T$, and fixed strike price 
$K$ is then the contingent claim on $[t_0, T]$  paying
$$ \left( {J(T)\over T\!-\!t_0}-K\right)^+
:=\max\left\{0,\, {J(T)\over T\!-\!t_0}-K
\right\}$$
at time $T$. Recall that points in time are taken to be non--negative real 
numbers. 
\goodbreak
\bigskip
{\bf 3.\quad The valuation problem for arithmetic--average Asian options:}\quad
Using the arbitrage--pricing principle, the price $C_{t, T}(K)$ at any time 
$t$ between $t_0$ and $T$ of the arithmetic--average Asian option introduced 
in paragraph two is given by:
$$C_{t, T}(K)=e^{-r(T-t)} E_t^Q\left[\left( 
{J(T)\over T-t_0}-K\right)^+\right].$$ 
From [{\bf GY}, \S 3.2] one has the following normalization 
of this value process:
\bigskip
{\bf Lemma:}\quad {\it For any time $t$ between $t_0$ and $T$, one has:
$$ C_{t, T}(K)={e^{-r(T-t)}\over T-t_0}\cdot {4S_t\over \sigma^2}\cdot
C^{(\nu)}(h,q), $$
where}
$$ C^{(\nu)}(h,q)=E^Q\biggl[ \max\bigg\{0,  
\int_0^h e^{2(B_\tau+\nu \tau)}d\tau
-q\bigg\}\biggr].$$
\smallskip
To explain the notation, the normalized parameters herein are as follows:
$$\leqalignno{
&\nu= {2\varpi\over \sigma^2}-1, &\cr 
h&={\sigma ^2\over 4}(T-t), &\cr 
q={\sigma^2\over 4S_t}&\left( K\cdot(T\!-\!t_0)
-\int_{t_0}^t S_\tau\, d\tau\right).
&\cr }$$
\smallskip
To interpret these quantities,  
$C^{\raise-1.5pt\hbox{$\scriptstyle (\nu)$}}(h,q)$ is the 
{\it normalized price\/} for the Asian option. It depends on the 
{\it normalized adjusted interest rate\/}  $\nu$, which for $\varpi$ 
positive is bigger than minus one, on the {\it normalized time to maturity\/}
$h$, which is non--negative, and on the {\it normalized strike price\/} $q$.   
\bigskip
To prove the Lemma break up the integral $J(T)$ at the point $t$ and bring 
$S_t$ in front of the expectation. Modify that part of the integral unknown 
given the information at time $t$ as follows. Change variables to obtain an 
integral from time $0$ to time $T\!-\!t$. The strong Markov property of 
Brownian motion identifies the exponent $B(t\!+\!\tau)-B(t)$ of the new 
integrand as a standard Brownian motion that is independent of the 
information %${\cal F}_t$ 
at time $t$. Herein change time by dividing through $\sigma^2/4$. Using the 
scaling property of Brownian motions, this time--changed process is again a 
standard Brownian motion if multiplied by $\sigma/2$ which completes the 
proof.    
\bigskip
{\bf 4.\quad Statement of results:}\quad In the setting of paragraphs one 
to three consider any time $t$ between the times $t_0$ and $T$. Using the 
decomposition
$$ C_{t, T}(K)={e^{-r(T-t)}\over T\!-\!t_0}\cdot {4S_t\over \sigma^2}\cdot
C^{(\nu)}(h,q)$$
of \S 3 Lemma, the valuation of the arithmetic--average Asian 
option is reduced to computing its normalized price $C^{(\nu)}(h,q)$. 
This price depends on the sign of the normalized strike price $q$. This
parameter is
recalled to be given by:
$$q
={\sigma^2\over4}\cdot  
{1\over S_{T-(4/\sigma^{2})\cdot h}}\cdot
\bigg( K\cdot (T\!-\!t_0)
-\int_{t_0}^{T-(4/\sigma^{2})\cdot h}
\hskip-2pt S_u\, du\bigg). 
$$
as a function of time $t$ which enters via the normalized time to 
maturity $h=(\sigma^2/4)\cdot (T\!-\!t)$. 
Generically $q$ will be positive. It is non--positive if 
and only if at time $t$ it is already known that the option will be 
in the money at its maturity $T$.  
\goodbreak
\bigskip
{\bf The main result:}\quad The main result of this paper is 
a closed form solution for $ C^{(\nu)}(h,q)$ if $q$ is positive. 
It expresses this function as the sum of integral representations.
These are obtained by integrating the product of Hermite functions, 
discussed in paragraph zero, with functions derived from weighted 
complementary error functions.   
Proved later in paragraphs nine to twelve,  the precise result is as 
follows: 
\bigskip
{\bf Theorem:}\quad {\it If  $q$ is positive, the normalized price 
$ C^{(\nu)}(h,q)$ of the Asian option is given by the following
five--term sum:} 
%\vskip3pt 
$$ %C^{(\nu)}(h,q)=
\ctrig +
{\Gamma(\nu\!+\!4)\over 2\pi(\nu\!+\!1)}\cdot
{(2q)\vbox to 9pt{}^{{\scriptstyle\nu+2\over\scriptstyle 2}}
\over
e\vbox to 9pt{}^{
{\scriptstyle 1\over\scriptstyle 2q}
+{\scriptstyle \nu^2h\over\scriptstyle 2}}
}
\cdot\Bigl( \chyp {\nu+2} +\chyp {-(\nu+2)}
-\chyp {\nu} - \chyp {-\nu} \Bigr).$$
\medskip
To explain the functions, recall from paragraph zero 
the Hermite functions $H_\mu$ of degree any complex number $\mu$. 
The {\it trigono\-me\-tric term\/} $\ctrig$ is then:
$$
\ctrig= c
\cdot
\int_0^\pi H_{-(\nu+4)} \left( -{\cos (\theta)\over \sqrt{2q}}\right)
\cos(\nu\cdot \theta)\, d\theta\, ,
$$  
recalling $\nu=2\varpi/\sigma^2\!-\!1$ and abbreviating:
$$
c=c(\nu,q)=
\Gamma(\nu\!+\!4)\cdot { e\vbox to 8pt{}^{2h(\nu+1)}-1\over \nu\!+\!1}
\cdot {(2q)^{{\scriptstyle\nu+2\over\scriptstyle 2}}
\over
\pi\cdot e^{\scriptstyle 1\over\scriptstyle 2q}}\, .
$$
For $\nu$ bigger than minus one $\ctrig$ has the additional representation:
\belowdisplayskip=0pt 
$$\ctrig
={ e\vbox to 8pt{}^{2h(\nu+1)}\!\!-\!1\over 2(\nu\!+\!1)}
\Bigl(1+2q(\nu\!+\!1)\Bigr)
+ 
\sin(\nu\pi)\cdot c\cdot \int\nolimits_0^\infty \hskip-4pt
H_{-(\nu+4)}\!\left({\cosh(x)\over \sqrt{2q\,} }\right)
e\vbox to 8pt{}^{-\nu x}\, dx \, .
$$
\goodbreak
\belowdisplayskip=8pt plus 5pt minus 3pt
The {\it hyperbolic terms\/} $\chyp b $ with 
$b$ equal to $\pm \nu$ or $\pm (\nu\!+\!2)$ are:
$$ \chyp b =
{2\over \sqrt{\pi}}  \cdot 
e\vbox to 9pt{}^{
{\scriptstyle \pi^2\over\scriptstyle 2h}+
{\scriptstyle b^2h\over\scriptstyle 2}}
\int_0^\infty 
H_{-(\nu+4)}\left( {\cosh (y)\over \sqrt{2q}}\right)
 \cdot e^{by} \cdot E_b(h)(y)\, dy\, . $$
Herein $E_b(h)$ is in terms of real--valued  functions for any real number 
$y$ given by: 
$$ E_b(h)(y)
%={e^{-2b^2h}\over \sqrt{2h\,}} \int_z^\infty 
%e^{-\left( {\scriptstyle u^2\over\scriptstyle 2h}+2ub\right) }
%\sin\left( \pi\left( 3b+{u\over h}\right)\right)du
=
\int_{ {\scriptstyle y\over \scriptstyle \sqrt{2h\,}}
+{\scriptstyle b\over\scriptstyle 2}\sqrt{2h\,}} ^\infty 
e^{-u^2}
\sin\biggl( \pi\biggl(b- u{\sqrt{2h\, }\over{ h\, }}\biggr)\biggr)du\, . 
$$
It is identified in paragraph ten as imaginary part of a certain weighted 
complementary error function at complex arguments.
Actually, this last interpretation seems to be best suited for 
a numerical analysis.
\goodbreak
\bigskip
{\bf The case $q\le 0$:}\quad It is a 
characteristic feature of Asian options that their prices become rather
rigid some time before maturity. The present case $q\le 0$ where at time 
$t$ it is known that the option will be in the money at its maturity
makes this precise. The rigidity of the price is reflected in the 
formula of [{\bf GY}, \S 3.4]. For any $h>0$ one has:
$$ C^{(\nu)}(h,q)={e^{2(\nu+1)h}-1\over 2(\nu+1)}-q\, .$$
%whence  
%$$ C_{t, T}(K)=S_t\cdot {1-e^{-r(T-t)}\over r(T-t_0)}-
%e^{-(T-t)}\left( K-{1\over T-t_0} \int_{t_0}^t S_\tau \, d\tau\right).
%\leqno{(10)}$$
The proof reduces to computing the $Q$--expectation of 
$\int_{\raise1.5pt\hbox{$\scriptstyle 0$}}^{
\raise-1.5pt\hbox{$\scriptstyle h$}}
\exp(2(B_\tau\!+\!\nu\tau))\, d\tau$. 
For an argument independent of
[{\bf GY}, \S 3.4], applying Fubini's theorem this expectation is the 
integral from zero to $h$ with respect to the variable $\tau$ of the 
product $\exp(2\nu\tau)$ times the expectation of $\exp(2B_\tau)$. This 
last expectation is the integral over the real
line with respect to the variable $x$ of the product $\exp(2x)$ times the 
density $(2\pi\tau)^{-1/2}\exp(-(2\tau)^{-1}x^2)$. Changing variables
$y=(2\tau)^{-1/2}(x-2\tau)$ it is seen to be equal to $\exp(2\tau)$. On 
substitution of this result one so is reduced to calculate the integral 
from zero to $h$ with respect  to $\tau$ of $\exp(2(\nu+1)\tau)$. 
Distinguishing the cases where $\nu+1$ is zero and 
non--zero, this is seen to have the above value, completing the proof. 
\goodbreak
\bigskip
{\bf Remark:}\quad The author is very grateful to a number of people, 
including a referee of this paper, for having drawn his attention to Yor's 
formula [{\bf Y}, (6.e), p.528]. 
His integral representation is a triple integral with a 
trigonometric--hyperbolic type integrand: 
$$C^{\raise -1.5pt\hbox{$\scriptstyle (\nu)$}}(h,q)
={ e\vbox to 10 pt{}^{
{\scriptstyle \pi^2\over\scriptstyle 2h}
-{\scriptstyle\nu^2 h\over\scriptstyle 2}}
\over \pi \sqrt{2h\, }}%\cdot
\int_0^\infty x^\nu\int_0^\infty 
 e\vbox to 10 pt{}^{-{\scriptstyle (1+x^2)y\over\scriptstyle 2}}
\cdot\Bigl( {1\over y}\!-\! q\Bigr)^{\hskip-2.pt +}
\cdot \psi_{xy}(h)\, dy\, dx\, ,$$
where for any positive real number $a$, the function $\psi_a$ is 
given for any $h>0$ by:
$$\psi_a(h)
=\int_0^\infty 
 e\vbox to 10 pt{}^{-{\scriptstyle w^2 \over\scriptstyle 2h}}
 e\vbox to 10 pt{}^{-a \cosh(w)}
\sinh(w)\cdot \sin\Bigl( {\pi\over h}w\Bigr) dw\, . $$
In contrast, our formula separates 
into a trigonometric and a hyperbolic part. It is given as a sum of
single integrals whose integrands have a structural interpretation as 
products of two functions. It moreover identifies the higher 
transcendental functions occuring as factors in these products, and  
shows how they are given by or built up from Hermite functions. 
%\medskip
On a technical level, these differences can be regarded as consequences 
of the different mathematical approaches for proving the valuation formula. 
%quite a characteristic feature for higher transcendental functions. 
%Due to the manifold mathematical approaches that are possible, it 
%actually is a characteristic feature of higher transcendental functions that 
%they have similar but different integral representations, and often 
%there are even quite a number of them. 
%\medskip
Which of the above representations is better suited for which purpose 
in which situation remains to be analyzed. 
\goodbreak
\vskip1.cm
\centerline{\gross Part II\qquad Proof of the valuation formula}
\vskip.35cm
{\bf 5.\quad  Computing the price of the Asian option using Laplace 
transforms:}\quad 
The basic idea for the Laplace transform approach to valuating Asian 
options is  not to consider a single above Asian option. Instead, take 
the exercise time $T$ of the option as variable and consider today at 
time $t$ all above Asian options whose exercise times $T$ range from $t$ 
to points in time very far in the future. Then, one wants to average 
the values of these options. This, however, is not possible in a naive 
way. One has to suitably weight these values and only then compute their 
average. This is effected by the Laplace transform.  
\medskip
In the generic case with $q$ positive a candidate for the Laplace 
transform of the factor $C^{(\nu)}(\ , q)$ of the price for the Asian 
option was computed in [{\bf GY}, pp.361ff]. Unfortunately, I found a 
conceptional error in their computations. Luckily however, Peter Carr 
eventually suceeded in making me see that this error is not as fatal 
as I originally believed. Indeed, in  his interpretation the 
results of [{\bf GY}] are sufficient for valuating Asian options 
as follows. 
\bigskip
{\bf Proposition:}\quad {\it  Suppose $\nu=2\varpi/\sigma^2-1>-1$.
If the positive real number $h$ is such that $q(h)=kh\!+\!q^*$ is 
positive, the normalized price of the Asian option at $h$ is given by:
$$C^{(\nu)}\Bigl(h,q\big(h \big)\Bigr)   %kh\!+\!q^*\bigr)
=\la^{-1}\Bigl( F_{GY}\big(q(h), z\big)\Bigr)\bigl(h\bigr)$$
as the Laplace inverse of $F_{GY}(q(h), z)$ at $h$.} 
\goodbreak
\bigskip
Leaving the function $q$ unexplained for a moment, for any positive 
real number $a$ define:
$$F_{GY}(a,z)={D_\nu(a,z)\over z\cdot\bigl(z-2(\nu+1) \bigr)}\, ,$$
for any complex number $z$ with positive real part bigger than 
$2(\nu\!+\!1)$, where on choosing the principal branch of the logarithm:
$$D_\nu(a,z)={e^{-{\scriptstyle 1\over\scriptstyle 2a}}\over a }
\int_0^\infty e^{-{\scriptstyle x^2 \over\scriptstyle 2a}} 
\cdot x^{\nu+3}\cdot  I_{\sqrt{ 2z+\nu^2}}\left( {x\over a}
\right)\, dx\, .$$    
Here $I_\mu$ is, for any complex number $\mu$,  the modified  Bessel 
function of order $\mu$ more fully discussed in paragraph six below or
in [{\bf L}, Chapter 5], for instance. 
%that for any
%complex number $z$ in ${\bf C}\1 {\bf R}_{<0}$ is given by the following 
%series:
%$$I_\mu(z)=\summe _{m=0}^\infty 
%{ %\left({\displaystyle z\over\displaystyle 2}\right)^{\mu+2m}
%1\over
% m!\, \Gamma(\mu+m+1)}\cdot 
%\left({\displaystyle z\over\displaystyle 2}\right)^{\mu+2m}
%\, ,$$
%where $\G(s)$ is the gamma function [{\bf L}, \S 5.7, (5.7.1)].   
\medskip
What regards the function $q$, the first general observation is that 
$T$ does not enter directly in the analysis but indirectly via the 
normalized time to maturity $h= (\sigma^2/4)\cdot (T\!-\!t)$ of paragraph 
three. Letting $T$ vary from 
$t$ to infinity, $h$ thus varies from zero to infinity. Take $h$ as new 
variable. The crucial difficulty now is that $q$ of paragraph three also 
depends on $h$ as a function of $T$. Indeed, abbreviating  $k=K/S_t$, on 
chasing definitions one gets:
$$ 
q(h)= k\cdot h + q^* \qquad\hbox{where}\qquad
q^*={\sigma^2\over 4S_t}\left( K\cdot(t\!-\!t_0)
-\int_{t_0}^t S_\tau\, d\tau\right).$$
Herein $q^*$ is independent of $T$ and thus of $h$. The upshot is that
$q$ is affine linear in $h$ and non--constant in particular. This is a 
nasty situation since the computations of the normalized value of the 
Asian option in [{\bf GY}, pp.361ff] erroneously take $q$ as a constant. 
Indeed, they compute for any positive real number $a$ and any complex 
number $z$ with sufficiently big positive real part the Laplace transform:
$$
F_{GY}(a,z):=\int _0^\infty e\vbox to 9pt{}^{-zx}
\cdot E\big[ A^{(\nu)}_x-a\big]\, dx\, ,
$$
setting: 
$$
A^{(\nu)}_x=\int_0^x e\vbox to 9pt{}^{ 2\cdot ( \nu \cdot u +B_u)}\, du\, , 
$$
%%{\raise 1pt\hbox{$\scriptstyle 0$}}^x\exp(2(\nu u +B_u))\, du $ 
for any non--negative real number $x$, to adopt Yor's notation.  
%However their computations
%are still sufficient for valuating Asian options as follows. 
%Indeed, in  the interpretation of Peter Carr these results give:
%
%\bigskip
%For the proof of the Proposition first define $F_{GY}(a,z)$ for 
%any complex number $z$ with sufficiently big positive real part 
%as the Laplace transform:
%$$
%F_{GY}(a,z)=\int _0^\infty e\vbox to 9pt{}^{-zx}
%\cdot E\big[ A^{(\nu)}_x-a\big]\, dx\, ,
%$$
%setting $A^{(\nu)}_x=\int_{\raise 1pt\hbox{$\scriptstyle 0$}}^x
%\exp(2(\nu u +B_u))\, du $ to accomodate Yor's notation.  In [{\bf GY}]
%these functions are proved to exist and to coincide with their 
%above definition. 
\bigskip
The proof of the Proposition then proceeds in two steps. Reduce 
to the above Geman--Yor Laplace transforms. Then sketch the essential 
steps of their computation.
\medskip
The basic idea of Peter Carr for reducing the Proposition to the 
Laplace transforms of [{\bf GY}] was to fix any positive $h$ and 
introduce an auxiliary function that coincides with the normalized 
value of the Asian option at $h$. More precisely, recalling \S 3 Lemma 
one has to consider the function $\psi$ given for any positive real number 
$x$ by:
$$\psi(x)=C^{(\nu)}\bigl(x,kx+q^*\bigr)
=E\big[ A^{(\nu)}_x-kx\!-\!q^*\big]\, .$$
Deform this function into the function $\phi$ given for any positive 
real number $x$ by:
$$\phi(x)=E\big[ A^{(\nu)}_x-kh\!-\!q^*\big]\, .$$
Then $\psi$ and $\phi $ are identical at $h$ by construction. 
Granted it exists, the Laplace transform of $\phi$ is equal 
to $F_{GY}(q(h),z )$ recalling $q(h)=kh\!+\!q^*$. Its Laplace inverse 
is the function $\phi$ by definition, and so equals $\psi$ when 
evaluated at $h$. One so has tautologically:
$$    
C^{(\nu)}\bigl(h,q(h)\bigr)=\psi(h)=\phi(h)
=\la^{-1}\bigl(F_{GY}(q(h),z )\bigr)(h)\, , $$
and is reduced to have the above Geman--Yor Laplace transforms 
$F_{GY}(a,z )$.  
\medskip
The basic idea for computing these Laplace transforms  in 
[{\bf GY}, pp.361ff] is to introduce a suitable stochastic clock for 
the process $J$ of paragraph two. This is made possible by the 
Bessel factorization result, attributed to 
Williams in [{\bf Y}, \S 1.5] and to Lamperti in [{\bf RY}, XI 1.28, p.432], 
that the exponential of any Brownian motion with drift 
is a time--changed Bessel process in the following way:
$$ e^{B_t+\nu t}=R^{(\nu)}\left( 
\int\nolimits_0^t e^{2(B_\tau+\nu \tau)}d\tau\right)
=R^{(\nu)}\big( A_t^{(\nu)}\big).$$
Here $R^{(\nu)}$ is the Bessel process on $[0,\infty)$ with index
$\nu$, starting at $0$ with the value $1$. At the level of a definition, the 
square of this process is  a continuous diffusion process $\rho$ with 
values in the non--negative real numbers satisfying:   
$$ d\rho_t=2(\nu+1)\,dt +2\sqrt{\rho_t}\, dB_t\, , \qquad \rho_0=1.$$
To be able to compare times, let $\tau_{\nu,u}$ % $\tau^{(\nu)}(u)$ 
for any non--negative real 
number $u$ denote the least upper bound of all non--negative real numbers 
$s$ such that 
$\int_{\raise1.5pt\hbox{$\scriptstyle 0$}}^{
\raise-.5pt\hbox{$\scriptstyle s$}}
 \exp(2(B_\tau\!+\!\nu \tau))\, d\tau =u$. Notice that
this equals 
$\int_{\raise1.5pt\hbox{$\scriptstyle 0$}}^{
\raise-.5pt\hbox{$\scriptstyle u$}}
 (R^{\raise-.5pt\hbox{$\scriptstyle (\nu)$}} (s))^{-2}\, ds$. 
\medskip
The computation of $F_{GY}(a,z)$ is then by 
deftly using the double role of $A^{(\nu)}$ as both the stochastic 
clock in its Bessel factorization and the actual underlying of the 
option. Indeed, fix any positive real 
number $x$,  and consider the process $A^{(\nu)}$ at $x$ on the set of
all events where the passage time $\tau_{\nu,a}$ %$\tau^{(\nu)}(a)$ 
takes values less than
or equal to $x$. Break the integral defing $A^{(\nu)}(x)$  at  
$\tau_{\nu,a}$. %$\tau^{(\nu)}(a)$. 
The first summand then is $A^{(\nu)}$ at time $\tau_{\nu,a}$ %$\tau^{(\nu)}(a)$
and so is equal to $a$. In the second summand  
one wants to restart the Brownian motion in the exponent of the integrand 
at $\tau_{\nu,a}$. %$\tau^{(\nu)}(a)$. 
Thus shift the variable of integration accordingly. 
The second integral then is the product of 
$\exp(2\cdot(B(\tau_{\nu,a}) \!+\!\nu\cdot \tau_{\nu,a}))$
%$\exp(2\cdot(B(\tau^{(\nu)}(a)) \!+\!\nu\cdot \tau^{(\nu)}(a)))$ 
times
%$A^{(\nu)}$ at $x\!-\!\tau^{(\nu)}(a)$, 
$A^{(\nu)}$ at $x\!-\!\tau_{\nu,a}$,
by abuse of language after having 
applied Strong Markov. This 
last process is such that it is independent of the information at time 
$\tau_{\nu,a}$. %$\tau^{(\nu)}(a)$. 
Recalling the role of $\tau_{\nu,a}$, %$\tau^{(\nu)}(a)$
the first factor 
is the square of the Bessel process $R^{(\nu)}$ at time $a$.  Now 
taking the expectation conditional on the information at $\tau_{\nu,a}$,
%$\tau^{(\nu)}(a)$,
one thus gets:
$$ E^Q\Bigl[  \bigl( A^{(\nu)}(x)-a\bigr)^+\, \Big|\,
 \fa _{\tau_{\nu,a}}\Bigr]
% \fa _{\tau^{(\nu)}(a)}\Bigr]
=\bigl( R^{(\nu)}(a)\bigr)^2\cdot  
E^Q\Bigl[A^{(\nu)} \Big(\big[x-\tau_{\nu,a}\big]^+\Bigr)\Big].$$
%E^Q\Bigl[A^{(\nu)} \Big(\big[x-\tau^{(\nu)}(a)\big]^+\Bigr)\Big].$$
%%
%%$$ E^Q_{\tau^{(\nu)}(q)}\biggl[  \Bigl( \int_0^he^{2(B_\tau+\nu \tau)} d\tau 
%%-q\Bigr)^+\biggr]=\bigl( R^{(\nu)}(q)\bigr)^2\cdot  
%%E^Q\biggl[\int_0^{h-\tau^{(\nu)}(q)}e^{2(B_\tau+\nu \tau)} d\tau
%%\biggr].$$
The  $Q$--expectation of 
$A^{(\nu)}(w)=
\int_{\raise1.5pt\hbox{$\scriptstyle 0$}}^{
\raise-.5pt\hbox{$\scriptstyle w$}}
\exp(2(B_u\!+\!\nu u))\, du$
is $(\exp(2(\nu\!+\!1)w)-1)/(2(\nu\!+\!1))$, as shown in \S 4 or 
using the more general results of [{\bf Y92}, pp.69ff]. It thus follows:
$$ E^Q\Big[ \big(A^{(\nu)}_x-a\bigr)^+\Big]
=E^Q\biggl[ \bigl( R^{(\nu)}(a)\bigr)^2\cdot  
{
e\vbox to 9pt{}^{\scriptstyle 2(\nu+1)\left[x-\tau_{\nu,a}\right]^+} -1
%e\vbox to 9pt{}^{\scriptstyle 2(\nu+1)\left[x-\tau^{(\nu)}(a)\right]^+} -1
\over  2(\nu+1)}\biggr].
$$
To avoid a direct calculation of this expectation adopt the common 
strategy to first calculate its Laplace transform at arguments $z$ whose 
real parts are positive and sufficiently big. The hope is to thus arrive at a 
simpler situation from which one is able to identify the original function
itself. There is a technical point in that one wants to interchange the 
expectation with the Laplace transform integral. If $z$ is real it is at 
this stage possibly best to follow Yor's proposal for justifying this. 
Indeed with the integrand of the double integral in question positive and 
measurable, apply Tonelli's theorem now but justify only in a later 
step that any of the resulting integrals is finite. The case of a general
argument $z$ is reduced to this case considering the absolute value of the 
integrand, and the result is:
$$ F_{GY}(a,z)
= {E^Q\bigl[ e^{-z\cdot \tau_{\nu,a}} \cdot\bigl( R^{(\nu)}(a)\bigr)^2
%= {E^Q\bigl[ e^{-z\cdot \tau^{(\nu)}(a)} \cdot\bigl( R^{(\nu)}(a)\bigr)^2
\bigr]\over z\cdot \bigl(z-2(\nu+1)\bigr)}\, .$$
One has enough information about Bessel processes 
for computing the $Q$--expectation of the right 
hand side. Using [{\bf Y80}, Th\'eor\`eme 4.7, p.80], 
the expectation of $ \exp(-z\cdot \tau_{\nu,a})$ %\tau^{(\nu)}(a))$ 
conditional upon 
$R^{(\nu)}(a)$ being the positive real number $w$ is given by the following
quotient of modified Bessel functions:
$$ E^Q\left[ e^{-z\cdot \tau_{\nu,a}%\tau^{(\nu)}(a)
}\Big\vert R^{(\nu)}(a)=w\right]
={I_{\sqrt{2z+\nu^2}}\over I_\nu}\left( {w\over a}\right).$$ 
This result ultimately depends  on the following standard fact about Bessel
processes [{\bf Y80}, (4.3), p.78]. If $\nu$ is bigger than minus one, 
the density $p_{\nu,a}(1,w)$ of the Bessel semigroup with index $\nu$ 
and starting point $1$ at time $a$ is given by:
$$ p_{\nu, a}(1,w)=
{w^{\nu+1}\over a}\cdot e^{-{\scriptstyle  1+w^2\over\scriptstyle 2a}}
\cdot I_\nu\left({w\over a}\right)
.$$
It is for this result that the hypothesis $\nu$ bigger
than minus one of the Proposition is crucially required. 
The numerator of the above Laplace transform now is obtained 
by integrating the product of these last two expressions 
and $w$ square with respect to  $w$ over the positive real line. With
the resulting integral finite for any complex number $z$ with real part 
bigger than $2(\nu\!+\!1)$, this completes the proof of the Proposition.
\goodbreak
\bigskip
{\bf 6.\quad Preliminaries on integral representations of gamma and Bessel
functions:} Recall the occurrence of modified Bessel functions in \S 5 
Proposition. They are  given for any complex number $\mu$ by the following 
series:
$$ I_\mu(z)=\summe _{m=0}^\infty  
{
%\left({\displaystyle z\over \displaystyle 2}\right)^{\mu+2m}
1\over m!\, \Gamma(\mu+m+1)} \cdot
\left({\displaystyle z\over \displaystyle 2}\right)^{\mu+2m} 
\, , $$
for any complex number $z$ in ${\bf C}\1 {\bf R}_{<0}$.
In this section preliminary results on  
integral representations of these functions are reviewed. 
\bigskip
{\bf Hankel contours and the Hankel representation of the gamma 
function:}\quad 
Recall that the gamma function at any complex number 
$s$ with $\re (s)>0$ is given by: 
$$\G(s)=\int_0^\infty e^{-x} x^{s-1}\, dx\, .$$
%for any complex number $s$ with $\re (s)>0$. 
In this section, following [{\bf D}, pp.225f], the Hankel form of 
the gamma function is described, an  integral representation of 
the reciprocal of $\G(s)$ valid for any complex number $s$. 
\medskip
First notice that the Laplace transform at 
any $z$ with $\re (z)>0 $ of the function $x^{s-1}$ on the 
positive real line is $z^{-s}\cdot\G(s)$. The following classical result
of Laplace now gives the inverse of this Laplace transform.
For any complex number $s$ with $\re (s)>0$ one has: 
$$ {1\over 2\pi i} 
\int_{\xi_0-i\infty}^{\xi_0+i\infty} e^{x\xi}\cdot \xi^{-s}\, d\xi
=
\la ^{ -1}\left(z^{-s}\right)(x)=
{x^{s-1}\over\Gamma (s)}\, ,$$
for any $x>0$ and any fixed $\xi_0>0$. Herein deform the path of integration
as follows. Fix any positive numbers $P$ bigger than $R$ and with $R$ bigger 
or equal $a$. Starting at $\xi_0-iP$ move parallel to the 
imaginary axis until $\xi_0+iP$. From that point move parallel to the 
real axis to the point $iP$ on the imaginary axis. Now follow counterclockwise
the circle with radius $P$ around the origin. When the parallel to the 
axis through $ia$ is hit follow this parallel until the circle with radius 
$R$ around the origin is hit. Continue clockwise on this circle until $R$
on the positive real axis.
Thereafter move to $\xi_0-iP$ in such a fashion that the path traced out in 
the complex plane is symmetric with respect to the real axis. 
\goodbreak
\medskip
The integral of $\exp(x\xi)\xi^{-s}$ over the closed path of integration just 
constructed is zero using the Cauchy Theorem. The integrand is rapidly 
decreasing in the radius $|\xi|$ of $\xi$. Letting $P$ go to infinity the
above inverse Laplace transform so equals the integral over the
{\it Hankel contour\/} $C_{a,R}$: the boundary of the pan in the complex plane
of radius $R$ around zero with  handle of diameter $2a$ stretching to minus 
infinity along the negative real line, 
\input xy
\xyoption{arc}
\xyoption{arrow}
$$
\xy     0;<1mm,0mm>                %%% setzt Koord--System: Einheit ist 1mm...
                \ar (10,-5);(30,-5),            %% (1)
                \ar @{-}(30,-5);(50,-5),        %% (2)
                \ar @{-*\dir{>}}(50,5);(30,5),  %% (3)
                \ar @{-}(30,5);(10,5),          %% (4)
                \ar @{.}(50,5);(55,5),	
                \ar 0;(72.06,0),                        %% die Ox -- Achse
                \ar (55,-12);(55,12),\ar@{},            %% die Oy -- Achse
                (55,5)*+\dir{*}*+!LU{\scriptstyle ia},
                (62.06,0)*+\dir{*}*+!UL{\scriptstyle R},
                (47.94,0)*+\dir{*}*+!UL{\scriptstyle -R},
                (55,0)*\cir<7.06mm>{dr^dl},\ar@{},

\endxy
$$
\vskip-.1cm
\centerline{\eightrm Figure A\quad 
The Hankel contour $\hbox{\eightit C}_{\scriptstyle a,R}$}
\goodbreak
\vskip.4cm
on which one comes in from $-\infty$ on the branch in 
$\{ z| \re(z), \im(z)<0\}$, passes counterclockwise around zero, and leaves 
on the branch in  $\{ z| \re(z)<0, \im(z)>0\}$. The Hankel contour $C_R$ is 
defined as the limit of $a$ going to zero of the Hankel contours $C_{a,R}$. 
\medskip
Returning to the gamma function, the main point now is that on any Hankel
contour $C$ the absolute value of $\exp(x\xi)$ is exponentially decreasing  
for any positive $x$ with the absolute value of $\xi$ going to infinity. 
Hence one obtains the following
{\it Hankel formula\/}:
$$ {x^{s-1}\over \Gamma(s)}
={1\over 2\pi i} \int _C e^{x\cdot\xi}\cdot\xi^{-s}\,  d\xi,$$
valid now for any complex number $s$, and any positive real number $x$. 
\bigskip
{\bf A Hankel type integral representation of Bessel functions:}\quad
In this section, following [{\bf WW}, $17\!\cdot\! 231$, p.362],  the 
following Hankel--type integral representations on the right half--plane 
of  modified Bessel functions are recalled. 
\bigskip
{\bf Lemma:}\quad {\it For any modified Bessel function $I_\mu$ one has:
$$  I_\mu(z)={1\over 2 \pi i }\int _{\log C}
e^{-\mu\cdot \xi+z\cdot\cosh(\xi)}\, d\xi, $$   
for any complex number $z$ with positive real part and any Hankel 
contour $C$.} 
\medskip
The key step in the proof is the following integral representation:
$$ I_\mu(z)={(z/2)^\mu
\over 2 \pi i } \int _C \xi^{-(\mu+1)}
\exp\Bigl({\xi+{ z^2\over 4 \xi}}\Bigr)
 \, d\xi, $$
valid for any complex number $z$ in ${\bf C}\setminus {\bf R}_{<0}$. 
Indeed, substitute  for the reciprocal gamma values in 
the series of  $I_\mu(z)$ the respective Hankel formulas 
with $x=1$. Interchange the order of summation and integration. The 
resulting factor $\summe_m (z^2/4\xi)^m/m!$ is the exponential function at 
$z^2/4\xi$, as was to be shown. 
\medskip
Change variables $\xi=z\eta/2$ in the above integral representation of 
$I_\mu$ to get:
$$  I_\mu(z)={1\over 2 \pi i } \int _C \eta^{-(\mu+1)}e^{
{\scriptstyle z\over 2}\left( \eta+{\scriptstyle 1\over
 \scriptstyle \eta}\right)} \, d\eta\, . $$
Herein change variables $\eta=\exp(\xi)$ to complete the proof of the Lemma.
\bigskip
{\bf The contour $\log C_R$:}\quad In the sequel,  it is the 
Hankel contours $C_R$ that are used. For a description of the 
logarithmicalized contour $\log C_R$, recall the principal branch of
the logarithm on the complex plane with the  non--positive real axis deleted:
$$ \log(z)=\log |z| + i\cdot \arg (z)
\qquad {\rm where}\quad -\pi <\arg (z)<\pi.$$
Returning to the Hankel contour $C_R$, the argument of the complex numbers 
in the upper branch of the panhandle in $C_R$ is $+\pi$, of those in the 
lower branch it is $-\pi$. The elements\goodbreak
on the circle part of $C_R$ 
have the form $R\cdot \exp(i\theta)$ with $-\pi<\theta<\pi$. The contour 
$\log C_R$ thus has the following shape: Coming in from plus infinity, move 
on the parallel through
$$
\xy     0;<1mm,0mm>            %%% setzt Koord--System: Einheit ist 1mm...
                \ar(0,0);(80,0),
                \ar(60,-10);(40,-10),
        \ar@{-}(40,-10);(20,-10)*+\dir{*}*+!RU{{\scriptstyle\log R-i\pi}},
        \ar@{-}(20,-10);(20,10)*+\dir{*}*+!RD{{\scriptstyle \log R+i\pi}},
                \ar(20,10);(40,10),
                \ar@{-}(40,10);(60,10),\ar@{},
                (20,0)*+\dir{*}*+!RU{{\scriptstyle\log R}}
\endxy
$$
\vskip-.25cm
\centerline{\eightrm Figure B\quad The contour $\scriptstyle\log C_R$}
\vskip.5cm
%Coming in from plus infinity, move on the parallel through 
the point $-i\pi$ to the real axis to the point $\log R-i\pi$. From this 
point move up to the point $\log R+i\pi$ on a parallel to the imaginary axis. 
Finally exit from  $\log R+i\pi$ to plus infinity on the parallel through 
$i\pi$ to the real axis. 
Notice $ \log C_R\subseteq \{ z| \re (z)>0\}$
if and only if the radius of the circle in $C_R$ is bigger than $1$. 
\goodbreak
\medskip
{\bf Odds and ends:}\quad To conclude, explicit 
consequences of the above developement are indicated. First of all  
using the explicit coordinates discussed above in the logarithmicalized 
Hankel contour with 
radius $R$ equal to one, the integral representation for $I_\mu$ 
of the Lemma specializes to the following 
{\it Schl\"afli integral representation\/}: 
$$ I_\mu(z)
= {1\over\pi}\int\nolimits _0^\pi e{}^{z\cos \theta} 
\cos (\mu\theta)\, d\theta
-{\sin(\mu\pi)\over \pi}\int\nolimits_0^\infty e{}^{-z\cosh(x)-\mu x}dx\, ,$$
for any complex number $z$ with positive real part [{\bf Wa}, $6\cdot 22$]. 
Direct calculations using the series of the respective modified Bessel 
functions or their Hankel type integral representations 
of the Lemma, prove the {\it recursion rule\/}:
$$ z\cdot I_\mu(z)=2(\mu+1)\cdot I_{\mu+1} (z)+z\cdot I_{\mu+2}(z)\, ,$$
valid for any complex numbers $z$ with $|\arg(z)|<\pi$ and 
$\mu$ [{\bf Wa}, $3\cdot 71$]. 
\goodbreak
\medskip
The  picture is finally completed with {\it Weber's integral\/}:
$$ \int\nolimits_0^\infty 
e{}^{-ax^2} x^{\mu+1} I_{\mu}(x)\, dx
= {1\over (2a)^{\mu+1}}\cdot e\vbox to 9pt{}^{
\scriptstyle 1\over \scriptstyle 4a},$$
valid for any positive number $a$ and any complex number $\mu$ with 
real part bigger than minus one. %$\re(\mu)>-1$. 
For its proof following [{\bf Wa}, $13\cdot 3$] develop the Bessel 
function factor of the integrand into its series. Interchange the 
order of integration and summation. Changing variables $y=ax^2$, any 
$n$--th term of the resulting series is
$(2a)^{-(\mu+1)}$ times the quotient of $(4a)^{-n}$ over $n$ factorial
times the quotient of 
$\int_{\raise1.5pt\hbox{$\scriptstyle 0$}}^{\raise-1pt
\hbox{$\scriptstyle \infty$}} \exp(-y)y^{\mu+n}\, dy$ 
over $\Gamma(\mu\!+\!n\!+\!1)$. With the real part of $\mu $ bigger than 
minus one, the numerator integral of the third quotient is the gamma 
function at $\mu\!+\!n\!+\!1$ and thus cancels with the denominator. 
Thus the series is $(2a)^{-(\mu+1)}$ times the series of the exponential 
function at the reciprocal of $4a$. Applying Lebesgue Dominated Convergence
then completes the calculation of Weber's integral.  
%It is proved by a direct calculation using the series of 
%the respective modified Bessel functions [{\bf Wa}, $13\cdot 3$]. 
\goodbreak
\bigskip
{\bf 7.\quad First steps of the Laplace inversion:}\quad  From paragraph 
five the Geman--Yor functions $F_{GY}$ are for any positive real number 
$a$ recalled to be given by: 
$$F_{GY}(a,z)
= {D_\nu(a,z)\over z\cdot\bigl(z-2(\nu+1) \bigr)}\, ,$$
where
$$D_\nu(a,z)={e^{-{\scriptstyle 1\over\scriptstyle 2a}}\over a }
\int_0^\infty e^{-{\scriptstyle x^2 \over\scriptstyle 2a}} 
\cdot x^{\nu+3}\cdot  I_{\sqrt{ 2z+\nu^2}}\left( {x\over a}
\right)\, dx\, ,$$
for any complex number $z$ with real part bigger than $2(\nu+1)$. 
In the sequel the problem of inverting this Laplace transform is
reduced to the following:
\bigskip
{\bf Lemma:}\quad {\it For $\nu$ bigger than minus one, the inverse Laplace 
transform at any positive real number $h$ of any function $F_{GY}(a,z)$ 
is given by:
$$
%
%\eqalign{
%\la& ^{\enspace -1}\left( {D_\nu(q,z)\over z(z-2(\nu+1))}  \right)(h)\cr  
%&= 
c_1\cdot\!\int_0^\infty \hskip-3pt{1\over 2\pi i}
\int_{ \log C_R }\hskip-13pt
e\vbox to 9pt{}^{
-{\scriptstyle x^2}
+{\scriptstyle x{\sqrt{\scriptstyle 2}\over\sqrt{\scriptstyle a}}
\cosh(w)} 
}
x^{\nu+3}\bigg\{\!
\la ^{ -1}\biggl(\! { e^{-w\sqrt{z\, }}\over z\!-\!(\nu\!+\!2)^2}\!\biggr)
-
\la ^{ -1}\biggl(\! {  e^{-w\sqrt{z\, }}\over z\!-\!\nu^2}  \biggr)
\!\bigg\}\!\left(\!{h\over 2}\!\right)  dw\,dx \, ,
%cr }
$$
where 
$$c_1=e\vbox to 9pt{}^{-h{\scriptstyle \nu^2\over\scriptstyle 2}}
\cdot 
{(2a)\vbox to 9pt{}^{\scriptstyle \nu+2\over\scriptstyle 2}
\over (\nu\!+\!1)\cdot
 e\vbox to 8pt{}^{\scriptstyle 1\over\scriptstyle 2a}}\, ,$$
where $C_R$ is any Hankel contour with $R\ge R_0$, and with integrands 
absolutely integrable.}
\goodbreak
\bigskip
The proof of the Lemma is based on the complex inversion formula 
for the Laplace transform reviewed in paragraph zero. 
\bigskip
To prove the Lemma, let  $z_0'$ be any positive real number such that 
the line $\{ z| \re (z)=z_0'\}$
is contained in the half--plane where $F_{GY}(a,z)$ is a 
holomorphic function. Suppose for a moment proved that this function 
is rapidly decreasing, and apply the complex inversion formula to it.  
Writing out $D_\nu(a,z)$, the inverse Laplace transform of $F_{GY}(a,z)$ 
is the function on the positive real line given for any positive real 
number $h$ by:
%$$
%\la ^{\enspace -1}\left( {D_\nu(q,z)\over z(z-2(\nu+1))}\right)(h)
%={1\over 2\pi i} \int _{z_0'-i\infty}^{z_0'+i\infty}
%e^{hz}{D_\nu(q,z)\over z(z-2(\nu+1))}\, dz,$$
% Substituting the expression for $D_\nu(q,z)$
%recalled above, this is equivalent to calculating the following integral:
$$ {e^{-{1\scriptstyle \over\scriptstyle 2a}}\over   a}\cdot 
{1\over 2\pi i} \int _{z_0'-i\infty}^{z_0'+i\infty}e^{hz}
\int_0^\infty 
{e\vbox to 9pt{}^{-{\scriptstyle x^2 \over\scriptstyle 2a}} 
\cdot x\vbox to 9pt{}^{\nu+3}\cdot  I_{\sqrt{ 2z+\nu^2}}\left( 
{\displaystyle x\over\displaystyle a}
\right)\over z(z-2(\nu\!+\!1))}
\, dx
\, dz \, .$$
For an equivalent expression, 
change variables $\eta =2z\!+\!\nu^2$, put $z_0=2z_0'\!+\!\nu^2$, 
and sub\-stitute the Hankel--type integral representation of \S 6 Lemma 
for the modified Bessel function. 
With $c_1'=4(\nu+1)\cdot (2a)^{-(\nu+4)/2}\cdot c_1$, 
the following integral is then to be computed:
$$c_1'\cdot
\biggl(\!{1\over 2\pi i}\!\biggr)^2\hskip-4pt
\int _{z_0-i\infty}^{z_0+i\infty} 
\hskip-4pt \int_0^\infty\hskip-4pt \int _{\log C_R}
\hskip-10pt
e\vbox to 9pt{}^{-{\scriptstyle x^2\over \scriptstyle 2a}}
x\vbox to 9pt{}^{\nu+3} 
e\vbox to 9pt{}^{{\scriptstyle x\over \scriptstyle a}\cosh(w)}\cdot 
e\vbox to 9pt{}^{{\scriptstyle h\over\scriptstyle 2}\cdot z}
{ e\vbox to 7pt{}^{\scriptstyle -w\sqrt{z}}\over (z-\nu^2)(z-(\nu\!+\!2)^2)}
\, dw\, dx\, dz\, . 
$$
\smallskip
%with $c_1'=4(\nu+1)\cdot (2q)^{-(\nu+4)/2}\cdot c_1$.
We claim that for $R$  sufficiently big the 
absolute va\-lue of the integrand of this integral is exponentially 
decreasing to zero with the absolute values of $z$, $x$, or $w$ going to 
infinity. 
\medskip
Granting this result, the above triple integral then gives the desired 
Laplace inverse, and, using Fubini's theorem, the order of its integrals 
can be interchanged. Take the integral for Laplace inversion, i.e., the 
integral over the line 
$\{ z| \re (z)=z_0\}$, as inner integral.  
Change variables $x=(2a)^{1/2}t$.  
The Lemma follows on decomposing the denominator of the integrand in partial 
fractions. 
%\goodbreak
\medskip 
One is thus reduced to prove the last claim about the asymptotic 
behaviour. For the calculations  
recall $|\exp(\xi)|=\exp(\re(\xi))$, for any complex number $\xi$. 
%\medskip
For the asymptotic behaviour in the absolute value of $w$, reduce to 
elements $w=x\pm i\pi$ in $\log C_R$. For these $\cosh(w)=-\cosh(x)$. Hence 
the absolute value of the hyperbolic cosine factor of the numerator
equals $\exp(-(x/a)\cosh(\re w ))$ whence the required asymptotic behaviour 
in the absolute value of $w$.
\medskip
For the asymptotic behaviour in $z$, notice 
$|\exp(-wz^{1/2})| =\exp(-\re(wz^{1/2}))$ and recall 
$z^{1/2}=\exp((1/2)\cdot(\log |z|+i\arg(z)))$. The argument of 
$z$ converges to $\pi/2$ with $|z|$ going to infinity. Now 
$\re(wz^{1/2})$ equals  $|z|^{1/2}$ times  
$\re (w)\cos(\arg(z)/2)- \im(w) \sin(\arg(z)/2)$. Herein the cosine is 
positive and bigger than $\cos(\pi/4)$. Thus 
$\re(wz^{1/2})$ is positive if the real part of any $w$ is big enough.  
Hence choose $R$ big enough to have the desired asymptotic behaviour
in the absolute value of $z$.
\medskip
The asymptotic behaviour in $x$ is  determined by $\exp(-x^2/2a)$, 
thus completing  the proof. 
% of the claim about the 
%asymptotic behaviour of the integrand in (7) and thus of the Lemma. 
%\medskip    
%the right hand side of (5) equals 
%$$
%\eqalign{
%C_1\cdot\int_0^\infty 
%e^{-{\scriptstyle x^2 \over\scriptstyle 2q}} 
%\cdot x^{\nu+3}\bigg\{
%{1\over 2\pi i} \int _{z_0-\infty}^{z_0+i\infty}
%e^{{\scriptstyle h\over\scriptstyle 2}\cdot z}
%{ I_{\sqrt{ z}}\left( {x\over q}
%\right)\over z-(\nu+2)^2}
%\, dz  
%-
%{1\over 2\pi i} \int _{z_0-\infty}^{z_0+i\infty}
%e^{{\scriptstyle h\over\scriptstyle 2}\cdot z}
%{ I_{\sqrt{ z}}\left( {x\over q}
%\right)\over z-\nu^2}
%\, dz \bigg\}\, dx\, , \leqno{(6)}$$
%with $C_1$ the constant of the Lemma, thus completing the proof. 
\goodbreak
\bigskip
{\bf 8.\quad Computation of certain Laplace transforms:}\quad Lacking a 
suitable reference, this section computes the inverse Laplace transforms 
identified in \S 7 Lemma.
\goodbreak
\medskip
For any $\alpha$ and $\beta$ in {\bf C}, consider the 
functions on the positive real line given by:
$$\leqalignno{
&f_{\alpha,\beta}(t)
={t^{-1/2}\over \sqrt{\pi\,}}
e\vbox to 9pt{}^{-{\scriptstyle\alpha^2\over\scriptstyle 4\cdot t}} 
-\beta\cdot 
e\vbox to 9pt{}^{\alpha\cdot\beta+\beta^2\cdot t}\, 
\er\!\left( \beta\sqrt{t\ }+{\alpha\over 2\sqrt{t\ }}\right),&\cr 
\noalign{\vskip4pt}
%\noalign{\goodbreak and}
g_{\alpha,\beta}(t)=& {e^{\beta^2\cdot t}\over 2}
\left( e^{+\alpha\cdot \beta}\, 
\er\!\left( {\alpha\over 2\sqrt{t\ }}+\beta\sqrt{t\ }\right)
+
e^{-\alpha\cdot \beta}\,
\er\!\left( {\alpha\over 2\sqrt{t\ }}-\beta\sqrt{t\ }\right)
\right),&\cr }
$$
for any positive real number $t$. 
%Suppose that $\alpha$ and $\beta$
%satisfy the following assumption:
%\bigskip
%\parein{40pt}{}
%{\it The real parts of the complex numbers $\alpha$ 
%and $\alpha^2$ are both positive. %,\hfil\break
%%and $\beta/\alpha $ is not a positive real number.
%}
%\paraus
%\bigskip
Then one has the following two results:
%\goodbreak
\bigskip
{\bf Lemma:}\quad {\it If the real parts of $\alpha$ and $\alpha^2$ 
are positive, one has:
$$ \la \left(f_{\alpha,\beta}\right)(z)
={e^{-\a\sqrt{z\ }}\over \sqrt{z\ }+\beta}\, ,$$
for any complex number $z$ in ${\bf C}\setminus{\bf R}_{<0}$ with real part 
bigger than $|\re (\beta)|$.}
\bigskip
{\bf Corollary:}\quad {\it If the real parts of $\alpha$ and $\alpha^2$ 
are positive, one has:
$$ \la \left(g_{\alpha,\beta}\right)(z)
={e^{-\a\sqrt{z\ }}\over {z\ }-\beta^2}\, ,$$
for any complex number $z$ in ${\bf C}\setminus{\bf R}_{<0}$ with real part 
bigger than $\re (\beta^2)$.}   
%\goodbreak
\bigskip
The Corollary follows from the Lemma upon decomposing the denominator 
in partial fractions
%$(z-\beta^2)^{-1}=(2\beta)^{-1}(1/(z^{1/2}-\beta)-1/(z^{1/2}+\beta))$
and using the linearity of the Laplace transform. 
\bigskip
We  use the following two results 
proved mutatis  mutandis in [{\bf D}, Beispiel 8, p.50f]:
$$\leqalignno{
\la \bigl( \psi_\alpha\bigr)(z)
=e^{-\alpha \sqrt{z\ }}
\qquad &\hbox{where}\qquad 
\psi_\alpha(t)=
{\alpha\over 2\sqrt{\pi\ }}\cdot  t^{-{3\over 2}}\cdot
e\vbox to 10pt{}^{-{\scriptstyle \alpha^2\over\scriptstyle 4\cdot t}}, 
&\cr 
\la  \bigl( \chi_\alpha\bigr)(z)={e^{-\alpha \sqrt{z\ }} \over \sqrt{z\ }}
\qquad &\hbox{where}\qquad 
\chi_\alpha(t)=
{1\over \sqrt{\pi\ }}\cdot t^{-{1\over 2}}\cdot
e\vbox to 10pt{}^{-{\scriptstyle \alpha^2\over\scriptstyle 4\cdot t}}, 
&\cr 
}$$
for any complex number $z$ in ${\bf C}\setminus {\bf R}_{\le 0}$ and 
any positive real number $t$. 
\medskip
Subtracting $\chi_\alpha$ from $f_{\alpha,\beta}$, the proof of the 
Lemma reduces to show the identity: 
$$
\la \left(
-\beta\cdot 
e^{\alpha\cdot\beta+\beta^2\cdot t}\, 
\er\!\left( \beta\sqrt{t\ }+{\alpha\over 2\sqrt{t\ }}\right)\right)(z) 
=
-\beta {e^{-\alpha \sqrt{z\, }}\over
\sqrt{z\,}(\sqrt{z\,}+\beta)}.$$ 
Multiplying any nice function $f$ with $\exp(a\cdot )$
induces a shift by $-a$ in its Laplace transform: 
$\la\, (\exp(at)f(t))(z)=\la\, (f)(z-a)$. Using this with $a=\beta^2$ one is 
further reduced to calculating the Laplace transform of the above
complementary error function factor only. 
\goodbreak
\medskip
Since $\re (\a)$  is positive, the real part of 
$\beta{u}^{1/2}+(\alpha/2){u}^{-1/2}$ goes to plus infinity 
with $u$ converging from the right to zero. 
The Fundamental Theorem of Calculus thus gives: 
$$\er\!\left( \beta\sqrt{t\ }+{\alpha\over 2\sqrt{t\ }}\right)
= -{2\over \sqrt{\pi\,}}
\int_0^t
e^{-\left( \beta \sqrt{u}+
{\scriptstyle \alpha\over \scriptstyle 2\sqrt{u}} \right)^2}
\left({\beta\over 2}u^{-1/2}-{\alpha\over 4}u^{-3/2}\right)\, du\, .
$$
\goodbreak
Using the transform--of--an--integral property 
$\la\, (\int_0^\bullet f(u)\, du)(w)=w^{-1}\la\,(f)(w)$ 
the La\-pla\-ce transform  of this complementary error function at 
$w=z-\beta^2$ is given by:
$$\eqalign{
& -{2\beta\over 2(z-\beta^2)} \la \left[ {1\over \sqrt{\pi\,}}\cdot
t^{-{1\over 2}}\cdot e^{-\left( \beta \sqrt{u}+{\scriptstyle \alpha
\over \scriptstyle 2\sqrt{u}}\right)^2}\,\right](z-\beta^2)\cr 
&+ {2\over 2(z-\beta^2)}  
\la \left[ {\alpha \over 2\sqrt{\pi\,}}\cdot t^{-{3\over 2}}
\cdot e^{-\left( \beta \sqrt{u}+{\scriptstyle \alpha
\over \scriptstyle 2\sqrt{u}}\right)^2}\,\right](z-\beta^2).\cr}
$$
Using the Laplace transform of $\chi_\alpha$, the first Laplace 
transform of this sum equals 
$${e^{-\alpha\beta}e^{-\alpha \sqrt{z\,}}\over \sqrt{z\,}}.$$ 
Using the Laplace transform of $\psi_\alpha$, the second Laplace 
transform of this sum  equals 
$$e^{-\alpha\beta}e^{-\alpha \sqrt{z\,}}.$$ 
The above sum hence equals  
$\exp(-\alpha\beta)\exp(-\alpha \sqrt{z\,})/(z^{1/2}(z^{1/2}+\beta))$. 
The identity to be proved follows upon substituting this last expression.
This completes the proof. 
\goodbreak
\bigskip
{\bf 9.\quad Two intermediate results:}\quad In this section the 
Laplace inversion of \S 7 is resumed concentrating on
the two summands of the expression derived~in~\S 7~Lemma.
\medskip
If $\nu\!>\!-1$, for any real numbers $a$, $h>0$ and $b$, thus define 
more generally:
$$ \ia^{ a}_{b}(h)
=
{1\over 2\pi i }\int_{\log C_R}
\int_0^\infty e^{-{\scriptstyle x^2}+ 
x {\sqrt{{\scriptstyle 2}\,}\over \sqrt{\scriptstyle a}\, } 
\cosh (w)}
x^{\nu+3}
\la ^{\enspace -1}\biggl( { e^{\scriptstyle - w\sqrt{z\,}}\over 
z-b^2    }     
\biggr)\left({h\over 2}\right)\, dx\, dw \,$$
with $\log C_R$ the logarithmicalized  Hankel contour discussed in \S 6. 
For these integrals  one has the following two results,  the first of a 
more structural and the second of a more technical nature:
\goodbreak
\bigskip
{\bf Lemma:}\quad {\it If $\nu\!>\!-1$, for any real numbers $R\ge 1$, 
and $a$, $h>0$, and $b$ one has:
$$
\ia^{ a}_{b}(h)
%\int_0^\infty e^{-{\scriptstyle x^2\over\scriptstyle 2q}} 
%x^{\nu+3}&
%\la ^{\enspace -1}\biggl( { I_{\sqrt{z\, }}
%\left({\displaystyle x\over\displaystyle q}\right)\over 
%$z-b^2    }     
%\biggr)\left({h\over 2}\right)\, dx\cr 
=
c_2\cdot\Gamma(\nu+4)\cdot {1\over 2\pi i}\int_{\log C_R} 
\Bigl( F_b(h)+F_{-b}(h)\Bigr)(w)
\cdot H_{-(\nu+4)} \Bigl( -{\cosh (w)\over \sqrt{2a\, }} \Bigr)\,  dw \, , 
%\cr
%\noalign{where}
%&C_2={ \Gamma(\nu+4)\cdot e^{{\scriptstyle h\over\scriptstyle  2}\cdot b^2 }
%\cdot (2q)^{\scriptstyle \nu+4\over\scriptstyle 2}\over 2} ,
%& \cr }
$$
where on both sides the respective integrands are absolutely integrable.}
\bigskip
\medskip
{\bf Corollary:}\quad {\it If $\nu\!>\!-1$, for any real numbers 
$\rho=\log R\ge 0$, and $a$, $h>0$, and $b$ one has:
$$
\ia^{ a}_{b}(h)={c_2\over \pi} \int_0^\infty \!
e^{-{\scriptstyle x^2}} 
x^{\nu+3}\bigg\{ \int_0^\pi \!\!
\re\bigl( G^a_{x,b}(h)\bigr)(\rho+i \theta )\, d\theta
+\!\int_\rho^\infty \!\!\im\bigl(G^a_{x,b}(h)\bigr)(y+i \pi )\,dy\bigg\}\, dx\,
 ,$$
with $\re(\xi)$ and $\im (\xi)$ the real respectively the imaginary 
part of any complex number $\xi$.}
\goodbreak
\bigskip
To explain the notation in these two results, $c_2$ denotes the following 
constant:
%Moreover recall from (4.1) 
%$ \er (z)=N(-z)=  (2/\sqrt{\pi})\int _z^{\infty} e^{-u^2}du
%=e^{-z^2}{2\over \sqrt{\pi}}\int_0^\infty e^{-u^2-2zu}du
%$, %\leqno{(1)}$$
%for any complex number $z$, and  
%the Hermite functions $H_\mu$ with $\re (\mu)<0$ that for any complex number 
%$z$ in ${\bf C}\1 {\bf R}_{<0}$ are given   by
%$H_\mu(z)={( \Gamma(-\mu))^{-1}}\int_0^\infty e^{-u^2-2zu}u^{-(\mu+1)}\, du
%$.
%\leqno{(1)}$$
%the functions $N$, $\er $ and $H_\mu$. Abbreviate: 
$$c_2=
c_2(h,b)={1\over 2}\cdot e^{h \cdot{\scriptstyle b^2\over\scriptstyle 2} }
. 
$$
The Hermite functions $H_\mu$ of degree $\mu$ are discussed in 
paragraph zero, and the two other functions that occur are given  by:
$$\displaylines{
F_b(h)(w)=e^{wb}\cdot 
\er \left( {w\over \sqrt{2h\,}}+{b\over 2}\sqrt{2h\,}\right),
\cr 
\noalign{\vskip3pt}
G^a_{x,b}(h)(w)=e\vbox to 9pt{}^{x
{\sqrt{{\scriptstyle 2}\,}\over \sqrt{\scriptstyle a}\, } 
\cosh (w)}
\Bigl( F_b(h)(w)+F_{-b}(h)(w)\Bigr),\cr 
%\noalign{and abbreviate:}
%&C_2={ e^{{\scriptstyle h\over\scriptstyle  2}\cdot b^2 }
%\cdot (2q)^{\scriptstyle \nu+4\over\scriptstyle 2}\over 2}. 
%&(4)\cr 
}
$$
for any complex number $w$.
\bigskip
{\bf Proof of the Lemma:}\quad For the proof of the Lemma  
choose a Hankel contour $C_R$ with $R$ so big that for any 
element $w$ in $\log C_R$ also the real parts of $w$ and $w^2$ are both 
positive. At any $w$ on $\log C_R $ substitute in 
$\ia^{ a}_{b}(h)$ 
the inverse Laplace transforms calculated in \S 8 Corollary.
Interchange the order of integration using the absolute integrability 
of the integrand. This gives 
the expression of the Lemma for $\ia^{ a}_{b}(h)$. The 
integrand herein is a holomorphic function on ${\bf C}\setminus {\bf R}_{<0}$.
Using the Cauchy Theorem, the value of the integral is  
independent 
of the Hankel contour $C_R$ chosen as long as $\log R$ is non--negative. This 
completes the proof.
\goodbreak
\bigskip
{\bf Proof of the Corollary:}\quad Given the Lemma, the proof of the 
Corollary is an exercise in path integration. Put $G_x=G^a_{x,b}(h)$. 
Changing the order of integration in the Lemma,
$$ c_2\int_0^\infty e^{-{\scriptstyle x^2}} 
x^{\nu+3} \, {1\over 2\pi i}\int_{\log C_R} G_x(w)\, dw\, dx$$
is to be calculated. Concentrate on its inner integral and consider the 
following subpath
$$ P=P_{c}+P_\infty \qquad {\rm with}\qquad P_{c}=P_{c,+}-P_{c,-}$$
of the path $\log C_R$. Here the path $P_{c,-}$ starts from $\log R$ and 
moves parallel to the imaginary axis to the point $\log R- i\pi$. The path 
$P_{c,+}$ starts from $\log R$ and moves parallel to the imaginary axis 
to the point  $\log R +i\pi$. The path $P_\infty$ moves from $\log R +i\pi$ 
parallel to the real axis to $+\infty$. For any $x>0$, the inner integral 
thus breaks up as follows:
$$
%\eqalign{
% &{1\over 2\pi i}\int_{\log C_R} G_x(w)\, dw\cr 
%&=
{1\over 2\pi i}\int_{P_\infty} \Bigl( G_x(w)-G_x(\ol w) \Bigr)\, dw    
+{1\over 2\pi i}\int_{P_{c,+}} G_x(w)\, dw\
-{1\over 2\pi i}\int_{P_{c,-}} G_x(w)\, dw\,
,
%\cr}
$$
with $\ol w $  the complex conjugate of $w$. Using
the series expansion of the exponential and the complementary error 
functions, $G_x$ is compatible with complex conjugation, i.e., $G_x$ 
evaluated at the complex conjugate of any complex number $w$ is the 
complex conjugate of $G_x$ at $w$:
$$ G_x\left( \ol w\right)=\ol{ G_x(w)}.$$
Write the elements of $P_\infty$ as $w=y+i \pi$ with $y\ge \log R$, 
notice $dw=dy$ and change variables. Using the compatibility
of $G_x$ with complex conjugation, 
one obtains  the improper integral from $\log R$ to 
infinity of the imaginary parts of $G_x(y+i \pi)$.
%
%$ in  it so follows:
%$$
%{1\over 2\pi i}\int_{P_\infty} \Bigl( G_x(w)-G_x(\ol w) \Bigr)\, dw 
%={1\over \pi} \int_{\log R}^\infty \im \left( G_x(y+i \pi)\right)\, dy\,.
%$$
\medskip
On the circle part $P_c$ of $P$ one has to be a bit more careful about the 
volume forms. The elements of $P_{c,+}$ are parametrized by 
$w=\log R +i \theta$ with $0\le \theta\le \pi$, whereas those of $P_{c, -}$ 
are parametrized by $\log R - i \theta $ with $0\le \theta\le \pi $. 
Changing variables accordingly and integrating from zero to $\pi$, the 
induced volume form in the $P_{c,+}$--integral is $dw=i\, d\theta$, whereas  
that on the $P_{C, -}$--part is $dw=- i\, d\theta$.
% Upon changing variables,the two parts of the integral over $P_c$ in (6) 
%thus are actually added.  
%As a consequence, 
%$${1\over 2\pi i}\int_{P_{c}} G_x(w)\, dw\
%={1\over \pi} \int_0^\pi \re\left(G_x(\log R +i \theta)\right)\, d\theta.
%$$
Using the compatibility of $G_x$ with complex conjugation, this completes 
the proof of the Corollary.
\goodbreak
\bigskip
{\bf 10.\quad The final explicit calculations:}\quad The aim of this 
section is to explicitly compute the integrals of \S 9 Corollary thus 
essentially proving the valuation formula \S 4 Theorem. 
\medskip
Recall the integral 
$\ia^{a}_{b}(h)$ of paragraph nine, the functions  $\er $, 
$H_\mu$ introduced in  paragraph zero,  and the function $E_b(h)$ 
introduced in  paragraph four. 
In $c_2(h,b)$ of paragraph nine divide by $\pi$, add a gamma factor, 
and drop the exponential in $h$ and $b$ to obtain:
$$ c_3= c_3(\nu)={1\over 2\pi}\cdot\Gamma(\nu+4)
%\cdot e^{{\scriptstyle h\over\scriptstyle 2}\cdot b^2 }
%\cdot (2q)^{\scriptstyle \nu+4\over\scriptstyle 2}
. $$
In view of \S 7 Lemma, the following result is the key for proving the 
valuation formula:
\goodbreak
\bigskip
{\bf Lemma:}\quad {\it If $\nu\!>\!-1$, for any real numbers $a$, $h\!>\!0$ 
and $b$, the function $\ia^a_{b}(h)$ is given by:}
$$\leqalignno{
&2\cdot c_3\cdot
e^{h{\scriptstyle b^2\over\scriptstyle 2}} \cdot
%\Big\{ N\Bigl({b\over 2}\sqrt{2h}\Bigr)
%    +N\Bigl(-{b\over 2}\sqrt{2h}\Bigr)
%\Big\}
%\cdot
\int_0^\pi H_{-(\nu+4)} \left( -{\cos (\theta)\over \sqrt{2a}}\right)
\cos(\nu\theta)\, d\theta
&\cr
\noalign{\vskip5pt} 
+ 
c_3\cdot{2\over \sqrt{\pi}} &\cdot e^{{\scriptstyle \pi^2\over\scriptstyle 2h}+
h {\scriptstyle b^2\over\scriptstyle 2}}\cdot
\int_0^\infty 
H_{-(\nu+4)}\left( {\cosh (y)\over \sqrt{2a}}\right)
\Bigl(e^{yb}E_{b}(h)(y)+e^{-yb}E_{-b}(h)(y)\Bigr) \, dy
\, .&\cr 
}$$
\goodbreak
\medskip
The Lemma is proved by computing the two integrals of $G_x:=G^a_{x,b}$ 
of \S 9 Corollary upon choosing $R=1$ there. Start with the following 
integral: 
$$ \!\int_0^\infty \!\!\im\bigl( G_{x}\bigr)(y+i \pi )\, dy\, . 
$$
Abbreviate $\beta_\pm= y/\sqrt{2h}\pm (b/2)\sqrt{2h\, }$ and notice 
$\cosh(y+i\pi)=-\cosh(y)$. Write the com\-ple\-mentary error functions 
occuring in $G_x(y+i\pi)$ as improper integrals starting from zero as in 
paragraph zero. Multiply out the expressions of the exponents to obtain:
$$ e^{\pm (y+i\pi)b}\, \er\Bigl( \beta_\pm+i{\pi\over \sqrt{2h\, }} \Bigr)
={2\over \sqrt{\pi\, }}\cdot e^{{\scriptstyle \pi^2\over\scriptstyle  2h}} 
\cdot e^{\pm yb}
\int_0^\infty e^{-(u+\beta_\pm)^2}\cdot 
e^{i\pi\left(\pm b -(u+\beta_\pm){\scriptstyle \sqrt{2h}\, 
\over\scriptstyle h }\right)}
du\, .$$
Since $\beta_\pm$ is a real number and $u$ can be taken as real numbers, 
the imaginary parts in $G_x$ are determined by the imaginary parts of 
the exponentials in the integral. These are given by the functions 
$E_{\pm b}(h)$ as indicated, completing the calculation.  
\goodbreak
\medskip
As a next step, 
$$\int_0^\pi \!\!
\re\bigl( G_{x}\bigr)(i \theta )\, d\theta  $$
is calculated. 
In contrast to the above argument, in the complementary error 
functions occuring in $G_{x}(i \theta ) $  now the following paths of 
integration are used: Abbreviating $\beta_\pm =\pm ({b/ 2})\sqrt{2h \,}$, 
first move from $\beta_\pm+i\theta $ to $\beta_\pm$, then continue from 
$\beta_\pm$ to plus infinity along the real line. 
Since $\cosh(i \theta)= \cos(\theta)$, the    
above integral then equals:
$$\leqalignno{ 
&\int_0^\pi e^{x {\sqrt{2} \over \sqrt{q\,}}\cdot \cos(\theta) } 
\Big\{
\bigl( \er (\beta_+)+\er (\beta_-)\bigr)\cdot \cos(\theta b)
+
{2\over \sqrt{\pi}}\bigl(%\phi_+(\theta)+\phi_-(\theta)
\phi_+ +\phi_-
\bigr)(\theta)\Big\}\, 
d\theta,  &\cr }  
$$
upon abbreviating for any angle $\theta$: 
$$
\phi_\pm(\theta) = -\re\biggl( e^{\pm i\theta b } 
\int _{\beta_\pm }^{\beta_\pm 
+
{\scriptstyle i\theta\over \scriptstyle \sqrt{2h\,} }} e^{-u^2}\, du\biggr).
 $$ 
In the first of these two last integrals notice 
$\er (\beta_+)+\er (\beta_-)=2$ since $\beta_-$ is minus $\beta_+$.
To calculate $\phi_\pm$ change variables $u=\beta_\pm+iw/\sqrt{2h\,}$ in the 
integral. Write the factor $i/\sqrt{2h\, }$ that is picked 
up as $\exp(i(\pi /2))/\sqrt{2h\,}$. Multiplying out the expression obtained 
in the exponent,  it follows
$$
e^{\pm i\theta b } \int _{\beta_\pm }^{\beta_\pm +
{\scriptstyle i\theta\over \scriptstyle \sqrt{2h\,} }} e^{-u^2}\, du\
= 
{e^{- {\scriptstyle b^2h\over \scriptstyle 2}}\over \sqrt{2h\,}}
\int_0^\theta e^{\scriptstyle u^2\over \scriptstyle {2h}}\cdot 
e^{ i\left( 
{\scriptstyle\pi\over\scriptstyle 2}
\pm b \left( \theta- u \right)
\right)}
du\, .$$
The real part of this expression is determined by the real part of
the exponential functions in the integral. 
Abbreviating $x_u=b(\theta-u)$, the values of the 
cosine at $\pi/2\pm x_u$ thus appear as factors. A shift by $\pi/2$ turns 
a cosine into a sine as follows:
$\cos(\pi/2\pm x_u)=\mp \sin (x_u)$. Thus $\phi_-$ is minus $\phi_+$,
%Hence  the second integral in  (4) is zero, 
completing the proof.
\goodbreak
\bigskip
{\bf 11.\quad First part of the proof of the valuation formula:}\quad 
The proof of the valuation formula of \S 4~Theorem is in two steps. As a 
first step,  \S 4~Theorem is in this paragraph established for $\nu$ 
bigger than minus one. As a second step, 
this equality is extended in the next paragraph to any complex number $\nu$ 
using analytic continuation. 
\medskip
Thus let $\nu$ be bigger than one. Recalling 
\S 7 Lemma, the valuation formula of \S4~Theorem then is obtained 
by subtracting $\ia^{ q}_{ \nu}(h)$ 
from $\ia^{ q}_{\nu+2}(h)$ 
and thereafter multiplying this difference with the constant $c_1$. 
Substitute the expressions computed in \S 10 Lemma.
Using $\exp(h(\nu+2)^2/2)=\exp(h\nu^2/2)\exp(2h(\nu+1))$ with the 
trigonometric summands one is reduced to show that
$$
\ctrig^*=\int_0^\pi H_{-(\nu+4)} \left( -{\cos (\theta)\over \sqrt{2q}}\right)
\cos(\nu\theta)\, d\theta\, 
$$
upon multiplication with  $c$ gives the first two terms in the 
sum of \S4~Theorem. 
As a first step, write $\ctrig^*$ as a double integral using the defining
integral representation of the Hermite function factor of its integrand. 
Applying Fubini's theorem and interchanging the order of integration 
the integrand of its inner integral is 
$\exp(2x(2q)^{-1/2}\cos \theta)\cdot \cos(\nu\theta)$. Substitute for this 
last integral using Schl\"afli's 
integral representation of paragraph six. Then reverse the order of 
integration to get:
$$ \ctrig^*
= {\pi\over \G(\nu+4)} \cdot I + \sin(\nu \pi)\cdot
\int\nolimits_ 0^\infty \hskip-4pt
H_{-(\nu+4)}\!\left( {\cosh(x)\over \sqrt{2q\ }}\right) e{}^{-\nu x}dx\, ,$$
where 
$$
I= 
\int\nolimits_0^\infty \hskip-4pt e^{-x^2} x^{\nu+3}
I_{\nu}\!\left(\! {2x\over \sqrt{2q\ }}\!\right)dx\, . 
%\qquad{\rm and}\qquad 
%%\cr  \noalign{\vskip4pt}
%J=\int\nolimits_ 0^\infty \hskip-4pt
%H_{-(\nu+4)}\!\left( {\cosh(x)\over \sqrt{2q\ }}\right) e{}^{-\nu x}dx\, .
$$
%Getting the asymptotic expansion of $\ctrig$ is thus reduced
%to calculating  $I$ and deriving the 
%asymptotic expansion of $J$. 
%\medskip   
The proof thus reduces to show:
$$
 I
={1\over 2}\cdot e\vbox to 9pt{}^{\scriptstyle 1\over\scriptstyle 2q}
\cdot (2q)\vbox to 9pt{}^{-{\scriptstyle\nu+2\over\scriptstyle 2}}\Bigl( 
1+2q\cdot(\nu+1)\Bigr).$$
Indeed, abbreviating $a=q/2$, change variables $w=a^{-1/2}x$ in $I$ to obtain: 
$$
I =
a\vbox to 9pt{}^{{\scriptstyle\nu+3\over\scriptstyle 2}
+{\scriptstyle 1\over\scriptstyle 2}}
\int\nolimits_0^\infty e^{-aw^2}w^{\nu+3}I_\nu(w)\, dw\, .$$
Apply the the recursion relation for modified Bessel functions 
of paragraph six:
$$I= 
a\vbox to 9pt{}^{\scriptstyle\nu+4\over\scriptstyle 2}\left[ 2(\nu+1)
\int\nolimits_0^\infty e^{-aw^2}w^{\nu+2}I_{\nu+1}(w)\, dw
+\int\nolimits_0^\infty e^{-aw^2}w^{\nu+3}I_{\nu+2}(w)\, dw\right]
. $$
The evaluation 
of $I$ is thus reduced to evaluating two Weber's integrals as reviewed in 
paragraph six. Substituting their respective values, the above identity 
follows. This completes the first step of the proof of \S 4 Theorem  
\goodbreak
\bigskip
{\bf 12.\quad Second part of the proof of the valuation formula:}\quad 
This second part of the proof of \S 4~Theorem extends its validity
from $\nu $ bigger than minus one, as established in the previous paragraph, 
to $\nu$ any complex number. 
\medskip
This reduces to show the following two results. First, the normalized 
price is an entire function in $\nu$. Second, the right hand side of 
\S 4~Theorem is a meromorphic function on the complex plane. Indeed, these 
two functions agree on real numbers $\nu$ bigger than minus one. Using 
the identity theorem they so agree 
on the complex plane as meromorphic functions. With one of them entire, 
the other one is entire, too. 
\medskip
In Yor's notation 
$A^{(\nu)}(h)=\int_{\raise1pt\hbox{$\scriptstyle 0$}}^{
\raise-1pt\hbox{$\scriptstyle h$}} 
\exp((2(B_u\!+\! \nu u))\, du$ recall 
$C^{\raise -1.5pt\hbox{$\scriptstyle (\nu)$}}(h,q)
=E[f(A^{\raise -2pt\hbox{$\scriptstyle (\nu)$}}(h))]$,  
where $f$ is given by $f(x)=(x\!-\!q)^+$, 
for any real number $x$. The proof of this normalized price being 
entire in $\nu$ is then further reduced to show 
$$
E\Big[f(A_h)\cdot e\vbox to 9pt{}^{\nu B_h}\Big]
$$ 
an entire function in $\nu$ where  $A$ is the 
process $A^{\raise -2pt\hbox{$\scriptstyle (0)$}}$. Indeed, this 
follows using the Girsanov identity
$E[f(A^{\raise -2pt\hbox{$\scriptstyle (\nu)$}}(h))]
=\exp(-\nu^2h/2)\cdot E[ f(A(h))\cdot\exp(\nu B(u))]$
%$$
%E\Big[f\Big(A^{\raise -2pt\hbox{$\scriptstyle (\nu)$}}_h\Big)\Big]
%=e\vbox to 9pt{}^{-\nu^2\cdot {\scriptstyle h\over \scriptstyle 2}}\cdot 
%E\Big[f\big(A_h\big)\cdot e\vbox to 9pt{}^{\nu B_h}\Big]
%$$
of [{\bf Y}, (1.c), p.510]. 
\goodbreak
\medskip
For proving the above function entire develop 
the factor $\exp(\nu B_h)$ of  $f(A_h)\cdot \exp(\nu B_h)$ 
into its exponential series. Suppose computing the expectation of the 
resulting series term by term is justified for any complex number $\nu$. 
For any complex number $\nu$, this then gives a power series in $\nu$ 
and thus explicitly shows $E[f(A_h)\cdot \exp(\nu B_h)]$ holomorphic at 
$\nu$.  
\medskip
Interchanging the order of integration and summation is 
justified using Lebesgue Dominated Convergence if the following is true. 
The expectations of the absolute value of any single term of the series 
for $f(A_h)\cdot \exp(\nu B_h)$ exist and the series  so obtained
converges. Using the Cauchy--Schwarz inequality this is implied by the series:
$$
 E\big[ f^2(A_h)\big]^{1/2}\cdot \summe_{n=0}^\infty {|\nu|^{n}\over n!}
\cdot E\big[ B^{2n}_h\big]^{1/2}
$$
being convergent for any complex number $\nu$. Implicit herein is 
that $f(A_h)$ is square integrable. This is implied by $A_h$ being
square integrable. Any $n$--th moment in particular of $A_h$ has 
been computed in [{\bf Y}, (4.d''), p.519] as: 
$$n! \cdot\bigg( {(-1)^n\over (n!)^2}+2\summe_{k=0}^n {(-1)^{n-k}\over
(n\!-\!k)!\cdot (n\!+\!k)!}
\cdot e\vbox to 9pt{}^{-h\cdot {\scriptstyle k^2\over \scriptstyle 2}}
\bigg).$$
In particular the second moment of $A_h$ is thus finite, as was to be 
shown.  What regards the even order moments of the Brownian 
motion $B$ at time $h$, they are computed as:
$$ E\big [ B_h^{2n}\big]= {(2h)^n\over n!}\cdot \Gamma\big(n+1/2\big),$$
for any non--negative integer $n$. The above series thus converges using
the ratio test. This completes the proof of the normalized price being an 
entire function in the paramater $\nu$. 
\medskip
The proof of the five--term sum of \S 4~Theorem being a meromorphic
function in  $\nu$ is based on the Hermite functions being
entire also in their degree. Indeed, from [{\bf L}, (10.2.8), p.285]
one has for any complex numbers $\mu$ and $z$ the representation:
$$H_\mu(z)={2^\mu\cdot\Gamma(1/2)\over \Gamma\big( (1\!-\!\mu)/2\big)}
\cdot\Phi\Big( -{\mu\over 2}, {1\over 2}, z^2\Big)
+
{2^\mu\cdot \Gamma(-1/2)\over \Gamma\big( -\mu/2\big)}
\cdot\Phi\Big( {1\!-\!\mu\over 2}, {3\over 2}, z^2\Big).
$$
Herein the reciprocal of the gamma function  is entire by 
construction, and the confluent hypergeometric function $\Phi$ is entire 
in its first and third variable.      
\medskip
There is a localization principle for proving entire a function on 
the complex plane. Indeed, one has to show analyticity at any fixed 
complex number. For this, one can restrict the function to 
any relatively compact or compact neighborhood of this fixed complex number.  
\medskip  
Apply this localization principle to the factor of the trigonometric term:
$$
A(\nu)=\int_0^\pi H_{-(\nu+4)} \left( -{\cos (\theta)\over \sqrt{2q}}\right)
\cos(\nu\theta)\, d\theta\, , 
$$
and restrict $\nu$ to belong to any suffiently small compact 
neighborhood $U$ of any point fixed in the complex plane.
The above integrand is analytic as a function in $\nu$ and smooth
as a function in $\theta$. Thus the absolute values of its derivatives 
with respect to $\nu$ of any order are bounded on the product of $U$ 
and the closed interval between zero and $\pi$. In this local situation, 
using a standard consequence of Lebesgue Dominated Convergence, 
differentation of $A$ with respect to $\nu$ is thus by partial 
differentiation with respect to the parameter $\nu$ under the integral 
sign. This proves $A$ entire as a function 
in $\nu$. The trigonometric term $\ctrig$ is obtained by multiplying
$A$ with $c$. Dividing off from $c$ the gamma factor 
$\Gamma(\nu\!+\!4)$,  it has a removable singularity in $\nu=-1$. Thus 
$\ctrig$ is meromorphic in $\nu$ with at most simple poles in the 
integers less than or equal to minus four. 
\medskip
What regards the hyperbolic terms, the claim is that as a function
of $\nu$ they can be extended as analytic functions from $\nu$ bigger
than minus one to to the whole complex plane. For this question it 
is sufficient to consider the function $B$ given by:
$$
B(\nu)=\int_0^\infty 
H_{-(\nu+4)}\left( {\cosh (y)\over \sqrt{2q}}\right)
 \cdot e^{b(\nu)y} \cdot E_{b(\nu)}(h)(y)\, dy\, , $$
where $b(\nu)$ equals $\pm\nu$ or $\pm(\nu\!+\!4)$ and with: 
%$E_{\xi}(h)$ is for any real number $y$ given by: 
$$ E_{\xi}(h)(y)
%={e^{-2b^2h}\over \sqrt{2h\,}} \int_z^\infty 
%e^{-\left( {\scriptstyle u^2\over\scriptstyle 2h}+2ub\right) }
%\sin\left( \pi\left( 3b+{u\over h}\right)\right)du
=
\int_{ {\scriptstyle y\over \scriptstyle \sqrt{2h\,}}
+{\scriptstyle \xi\over\scriptstyle 2}\sqrt{2h\,}} ^\infty 
e^{-u^2}
\sin\biggl( \pi\biggl(\xi- u{\sqrt{2h\, }\over{ h\, }}\biggr)\biggr)du\, . 
$$
for any real number $y$. Extension of this function is essentially by
reduction to the case where the real part of the degree $-(\nu\!+\!4)$
of the Hermite function factor in the integrand of $B(\nu)$ is negative, 
or equivalently, the real part of $\nu$ is bigger than minus four.   
\medskip
To fix ideas first consider the case where $\nu$ is such that the
degree of the Hermite function factor in the integrand of $B(\nu)$ 
is a non--negative integer. This Hermite function is then the 
corresponding Hermite polynomial. In the integrand of 
$B(\nu)$ the absolute values of the Hermite function 
factor and the exponential function factor so have linear exponential 
order in the variable $y$. The decay to zero of the absolute value of 
the function $E_{b(\nu)}$, however, is  of square exponential order 
in $y$. It thus dominates the asymptotic behaviour with $y$ to infinity. 
The absolute value of the integrand of $B (\nu)$ 
is so majorized by an integrable function,  and  $B$ can 
be extended to the above values of $\nu$.  
\medskip
For the general case of the reduction, fix any $\nu$ of real part less 
than or equal to minus one that is not an integer. Apply the above 
localization principle and let $\nu$ belong to a sufficiently small 
compact neigborhood $U$ in the half--plane $\{ \re (z)\!\le \! -1\}$. 
Shrinking $U$ if necessary assume that for any $\nu$ in $U$ the degree 
of the Hermite function in the integrand of $B(\nu)$ is not an integer. 
Using the recursion rule for Hermite functions of paragraph zero express 
the Hermite function factor of $B(\nu)$ in terms of weighted Hermite 
functions of negative degrees. Further shrinking $U$ if necessary, assume 
that the so obtained relation represents $B$ on $U$. Herein the 
weighting factor for the respective Hermite functions are given as 
powers of $z=(2q)^{-1/2}\cosh(y)$  times polynomials in $\nu$. 
The absolute values 
of the  polynomials in $\nu$ can be majoriozed uniformly on $U$. The 
problem thus reduces to majorize by an integrable function on the 
positive real line in the variable $y$ finitely many functions on $U$ 
times the positive real line sending $\nu$ and $y$ to:
$$
H_{-(\nu+4+k)}\left( {\cosh (y)\over \sqrt{2q}}\right)
 \cdot\cosh^\l (y) \cdot e^{b(\nu)y} \cdot E_{b(\nu)}(h)(y)\, , 
$$ 
where $k$, $\l$ range over finitely many non--negative integers 
and $k$ is such that $\nu\!+\!4\!+\!k$ is positive. The leading terms
of the  asymptotic expansion of paragraph zero for any Hermite function 
$H_\mu(z)$ with degree $\mu$ any complex number with negative real part 
has order $z^\mu$. Asymptotically with $y$ to infinity, the Hermite 
function with the smallest positive number $\nu\!+\!4\!+\!k$ thus
dominates the other Hermite function factors in the above functions.
With $\nu$ ranging over a compact set, there is a minimal such degree
on $U$. Similarly, there are such majorizing choices $\l^*$ for  
the  factors $\cosh^\l (y)$, and $\nu^*$ for the absolute values of the 
factors $ \exp(b(\nu)y)$, and  $\nu^{**}$ for the absolute values of the 
factors $E_{b(\nu)}(h)(y)$. A  four--factor--majorizing function on the 
positive real line thus results whose asymptotic behaviour with $y$ to 
infinity is governed by the square--exponential decay to zero of the 
corresponding factor $E_{b(\nu^{**})}(h)$ and which is integrable.  
\medskip
If the real part of $\nu$ is bigger than minus one, the above argument
holds in a simplified form. The upshot so is that any complex number not 
an integer has a sufficiently small compact neighborhood such that the 
absolute value of the integrand of $B(\nu)$ on $U$ times the positive 
real line can be majorized by an integrable function on the positive real 
line. Herein, compact neigborhoods can be replaced by relatively compact
neigborhoods mutatis mutandis. Using Lebesgue Dominated Convergence, $B$ 
can thus be extended as a continuous function to the whole complex plane
with the integers less than or equal to minus four deleted.  
\medskip
The idea for showing $B$ analytic as a function of $\nu$ on the complex 
plane with the integers less than or equal to minus four deleted is as 
follows. Show that differentiation of $B$ is by differentation under 
the integral sign and use that its integrand is entire as function 
of $\nu$. For this again first localize to $\nu$ in any sufficiently 
small compact neighborhood containing no integers less than or equal to 
minus four. The aim is then to majorize the absolute value of the 
derivative with respect to $\nu$ of the integrand of $B$ by an integrable 
function independent of $\nu$ as above. The above argument for getting 
such a majorizing function is based on a comparison of decay rates. The 
integrand of $B$ has one factor which on the positive real line decays 
to zero of square exponential order whereas the other factors explode 
of at most linear exponential order. This situation is preserved on 
differentiation with respect to the parameter $\nu$. In particular, 
differentiating with respect to the degree the asymptotic expansion 
for Hermite functions on the right half--plane gives an asymptotic 
expansion for this function's partial derivative with respect to the 
degree. 
\medskip
At this stage it remains to extend $B$ analytically to the intgers 
less than or equal to minus four. However, $B$ remains bounded in any 
punctured compact neigbourhood of such an integer. Thus $B$ can be 
extended to an entire function, completing the proof of \S 4~Theorem.
\goodbreak
%\bigskip
\vskip.8cm
\centerline{\gross Part III}
\vskip.3cm
{\bf 13.\quad Remarks about hedging:}\quad This paragraph's aim
is to compute the Asian option's Delta and discuss how the 
seller's hedging portofolio is determined by~it. 
\medskip
Delta,  and similarly 
the other local hedging parameters, are computed by partially 
differentiating the price function $C_{t,T}(K)$ of paragraph three:
$$ C_{t,T}(K)= 
e^{-r(T-t)}\cdot S_t\cdot {4\over \sigma^2(T-t_0)}\cdot
C^{(\nu)}(h,q). \leqno{} $$
To simplify notation, in the sequel as many arguments of a function as 
possible are suppressed, thus writing $C_t=C_{t,T}(K)$ and $C^{(\nu)}= 
C^{(\nu)}(h,q)$ in particular. 
\goodbreak
\medskip
Hedging of the Asian option as a particular case of the general 
theory of hedging European--style contingent claims in a complete 
Black--Scholes economy has been discussed in [{\bf K}, p.23f]. The 
seller's hedging portofolio $\Pi$ is determined using the martingale 
representation of the conditional expectation of the discounted pay--out 
of the Asian option at its time of maturity, i.e., using the following
stochastic differential equation:
$$ e^{-r(t-t_0)} C_{t,T}(K)= C_{t_0,T}(K)+
\sigma \int_{t_0}^t e^{-r(s-t_0)} \Pi_s \, dB_s\,. \leqno{}$$
Apply to the martingale of the left hand side the It\^ o formula.  
Comparing diffusion coefficients, the hedging portofolio $\Pi_t$ at any 
time $t$ between $t_0$ and $T$ thus is given by:
$$ \Pi_t=S_t\cdot \Delta_t
\qquad \hbox{where} \qquad 
\Delta_t={\partial C_{t}\over \partial S_t}.\leqno{}$$
The option's Delta $\Delta_t$ is computed in the case where $q$ 
is positive and  \S 4 Theorem applies. 
\medskip
The partial derivative of 
$q$ with respect to $S_t$ times $S_t$ being minus $q$, one has:
$$
{\partial C_{t}\over \partial S_t}
=e^{-r(T-t)}{4\over \sigma^2(T-t_0)}\Bigl( 
C^{(\nu)}(h,q)- q\cdot {\partial C^{(\nu)}\over \partial q}(h,q)\Bigr).
\leqno{}$$
One is thus reduced to computing the partial derivative with respect to  
$q$ of $C^{\raise -1.5pt\hbox{$\scriptstyle (\nu)$}}$ at $(h,q)$. Think 
of this as the option's normalized Delta. Further reduce as follows. Write 
the five--term sum of \S 4~Theorem for  
$C^{\raise -1.5pt\hbox{$\scriptstyle (\nu)$}}$ in the form:
$$C^{(\nu)}=d\cdot D_{\rm trig}+d\cdot D_{\rm hyp},$$ 
setting:
$$d=d(\nu,q)={\Gamma(\nu\!+\!4)\over 2\pi(\nu\!+\!1)}
\cdot { (2q)\vbox to 8pt{}^{ {\scriptstyle \nu+2\over \scriptstyle 2}}
\over e\vbox to 8pt{}^{ {\scriptstyle 1\over \scriptstyle 2q}}}
\, , $$
and where the modified trigonometric term $D_{\rm trig}$ and the 
modified hyperbolic term:
$$ e\vbox to 9pt{}^{{\scriptstyle \nu^2h\over \scriptstyle 2}}\cdot
D_{\rm hyp}
=\chyp {\nu+2} +\chyp {-(\nu+2)}
-\chyp {\nu} - \chyp {-\nu} $$
are functions in the variables $\nu$, $h$ and in $q=q(h)$. 
The partial derivative of $C^{\raise -1.5pt\hbox{$\scriptstyle (\nu)$}}$
with respect to $q$ then is:
$$    
%{\partial C^{(\nu)}\over \partial q}(h,q)= 
\Bigl( D_{\rm trig}+D_{\rm hyp}\Bigr)\cdot
{\partial d\over \partial q}+ 
d\cdot \Bigl({\partial D_{\rm trig}\over \partial q}
+{\partial D_{\rm hyp}\over \partial q}\Bigr).\leqno{}$$
Thus, one is further reduced to compute the partial derivatives with 
respect to $q$ of the above functions $d$, $D_{\rm trig}$, and $\chyp b $
with $b$ equal to $\pm \nu$ or $\pm(\nu\!+\!2)$. The first of these is:
$${\partial d \over \partial q}
= {\Gamma(\nu+4)\over \pi(\nu+1)}\cdot \Bigl(2+(\nu+2)2q\Bigr)\cdot
(2q)^{{\scriptstyle \nu-2\over \scriptstyle 2}}\cdot
 e^{-\,{\scriptstyle 1\over \scriptstyle 2q}}\,. \leqno{}$$
To compute the partial derivatives of the other two functions 
recall from paragraph zero the recursion formula for the derivative of 
Hermite functions.
Granting for a moment that one can justify the differentiation under the 
integral sign, the partial derivative with respect to  $q$ of 
$D_{\rm trig}$ at $(\nu,h,q)$ is: 
$$  
2\cdot\Bigl(1-  e^{2h(\nu+1)}
\Bigr) 
\cdot{2(\nu+4)\over (2q)^{3/2}}
\int_0^\pi H_{-(\nu+5)} \left( -{\cos (\theta)\over \sqrt{2q}}\right)
\cos(\nu\theta)\cdot \cos(\theta)\,  d\theta
.\leqno{}
$$ 
With the functions $E_b(h)$ of paragraph four independent of $q$,
the partial derivative with respect to the variable $q$ of any
$\chyp b $ at $(\nu,h,q)$ similarly is:
$$
{2\over \sqrt{\pi\, }}
e\vbox to 9pt{}^{\scriptstyle \pi^2\over\scriptstyle 2h}
\cdot 
{2(\nu\!+\!4)\over (2q)^{3/2}}%\cdot
\int_0^\infty 
H_{-(\nu+5)}\!\left(\! {\cosh (y)\over \sqrt{2q}}\right)\cdot\cosh (y)\cdot
e\vbox to 9pt{}^{yb}E_b(h)(y) \, dy
.\leqno{}$$
To justify the differentiations under the integral sign first notice 
that differentation is a local concept and thus the parameters $\nu$,
$h$ can be restricted to vary in a fixed compact set not containing
points with $q$ equal to zero. The integrands of $D_{\rm trig}$ and 
$D_{\rm hyp}$ then are integrable and differentiable functions not 
only in the variable $\theta$ respectively $y$ but also in the 
variables $\nu$, $h$, and in $q$. 
Consider the maxima in $h$ of the absolute value of the partial 
derivatives with respect to $q$ of $D_{\rm trig}$ and $D_{\rm hyp}$, 
for any fixed triple  of the other variables. In particular with $q$ 
bounded and bounded away from zero the so obtained  functions in the 
variables $\theta$ respectively $y$, $\nu$ then are majorized by 
integrable functions. A standard application of Lebesgue Dominated 
Convergence thus completes the argument for  the differentiations under 
the integral sign, and thus completes the calculations.  
\goodbreak
\medskip
\bigskip
{\bf References}
\medskip
\vbox{\baselineskip=9.5pt
\ninerm
\parein{40pt}{\ninebf [B]} R. Beals: {\nineit Advanced mathematical 
analysis\/}, GTM 12, Springer 1973
\paraus
\parein{40pt}{\ninebf [D]} G. Doetsch: {\nineit Handbuch der Laplace 
Transformation\/} I, Birkh\"auser Verlag 1971
\paraus
\parein{40pt}{\ninebf [DSM]} D. Duffie: {\nineit Security markets\/},
Academic Press 1988
\paraus
\parein{40pt}{\ninebf [GY]} H. Geman, M. Yor: Bessel processes, Asian  
options, and perpetuities, {\nineit Math. Finance\/} {\ninebf 3}(1993), 
349-375
\paraus
\parein{40pt}{\ninebf [K]} I. Karatzas: {\nineit Lectures on the mathematics
of finance}, CRM Monographs 8, American Mathematical Society, Providence 1997
\paraus 
\parein{40pt}{\ninebf [L]} N.N. Lebedev:  {\nineit Special functions and
their applications\/}, Dover Publications 1972
\paraus
\parein{40pt}{\ninebf [RY]} D. Revuz, M. Yor: {\nineit Continuous martingales 
and Brownian motion}, 2nd ed., Springer 1994
\paraus
\parein{40pt}{\ninebf [Sch]} M. Schr\"oder: On the valuation of 
arithmetic--average Asian options: explicit formulas, Universit\"at 
Mannheim, M\"arz 1999
\paraus
\parein{40pt}{\ninebf [WW]} E.T. Whittacker, G.N. Watson: {\nineit A 
course in modern analysis\/}, Cambridge UP, repr. 1965 
\paraus
\parein{40pt}{\ninebf [Y80]}  M. Yor: Loi d'indice du lacet Brownien, et
distribution de Hartman--Watson, {\nineit Z. Wahr\-schein\-lichkeitstheorie\/}
{\ninebf 53}(1980), 71--95
\paraus
\parein{40pt}{\ninebf [Y]}  M. Yor: On some exponential functionals of 
Brownian motion, {\nineit Adv. Appl. Prob.} {\ninebf 24}(1992), 509--531
\paraus   
\parein{40pt}{\ninebf [Y92]} M. Yor: {\nineit Some aspects of Brownian 
motion\/} I, Birkh\"auser 1992
\paraus
}
\vfill
\eject
\headline={  }
\hoffset-.3cm
\nopagenumbers
\baselineskip20pt
\centerline{  }
\vskip4truecm
\centerline{\gr  On the valuation of arithmetic--average }
\vskip2pt
\centerline{\gr Asian options: integral representations} 
\vskip.3cm
\baselineskip14pt
\centerline{\grossrm by}
\centerline{\gross M. Schr\"oder}
\centerline{\grossrm (Mannheim/Bonn)}
\vskip.3cm
\centerline{\grossrm October 1997}
\centerline{\grossrm (revised November 1999)}
  \bye

\bye